\documentclass[10pt,twoside]{article}

\usepackage[paperwidth=155mm,paperheight=235mm,
  left=19mm,right=19mm,top=21mm,bottom=16mm,
  headheight=11pt,headsep=13pt]{geometry}
\usepackage[T1]{fontenc}
\usepackage{amsmath,amsthm,mathtools}
\newcommand{\FockStandardFonts}{%
  \IfFileExists{mathptmx.sty}{\usepackage{mathptmx}}{}%
  \usepackage{amssymb}%
}
\IfFileExists{newtxtext.sty}{%
  \IfFileExists{newtxmath.sty}{%
    \usepackage{newtxtext,newtxmath}%
  }{\FockStandardFonts}%
}{\FockStandardFonts}
\usepackage[protrusion=true]{microtype}
\usepackage{enumitem,fancyhdr,booktabs,array,needspace,etoolbox}
\usepackage[colorlinks=true,linkcolor=blue,citecolor=blue,urlcolor=blue]{hyperref}
\hypersetup{
  pdftitle={Carleson Measures and Integral Operators on Fock-Sobolev Spaces},
  pdfauthor={Jianxiang Dong and Shibo Hu},
  pdfsubject={Critical symbols, essential norms, Fredholm indices, and differential equations on Fock-Sobolev spaces},
  pdfdisplaydoctitle=true}

\setlist[enumerate]{label=(\alph*),leftmargin=1.8em,itemsep=2pt,topsep=4pt}

\makeatletter
\def\@seccntformat#1{\csname the#1\endcsname.\enspace}
\renewcommand\section{\@startsection{section}{1}{\z@}{-3.1ex plus -1ex minus -.2ex}{1.2ex plus .2ex}{\normalfont\normalsize\bfseries}}
\renewcommand\subsection{\@startsection{subsection}{2}{\z@}{-2.2ex plus -.8ex minus -.2ex}{.8ex plus .2ex}{\normalfont\normalsize\bfseries}}
\makeatother

\allowdisplaybreaks[1]

\newtheoremstyle{journal}{7pt}{7pt}{\itshape}{}{\bfseries}{.}{.5em}{}
\theoremstyle{journal}
\newtheorem{theorem}{Theorem}[section]
\newtheorem{lemma}[theorem]{Lemma}
\newtheorem{proposition}[theorem]{Proposition}
\newtheorem{corollary}[theorem]{Corollary}
\theoremstyle{definition}

\newtheorem{example}[theorem]{Example}
\BeforeBeginEnvironment{theorem}{\Needspace{5\baselineskip}}
\BeforeBeginEnvironment{lemma}{\Needspace{5\baselineskip}}
\BeforeBeginEnvironment{proposition}{\Needspace{5\baselineskip}}
\BeforeBeginEnvironment{corollary}{\Needspace{5\baselineskip}}
\BeforeBeginEnvironment{example}{\Needspace{5\baselineskip}}
\BeforeBeginEnvironment{proof}{\Needspace{4\baselineskip}}

\fancypagestyle{firstpage}{%
  \fancyhf{}%
  }

\newcommand{\C}{\mathbb C}
\newcommand{\N}{\mathbb N}
\newcommand{\Nzero}{\mathbb N_0}
\newcommand{\F}{\mathcal F}
\newcommand{\Hc}{H(\C)}
\newcommand{\du}{\mathcal D_{\mathbf u}}

\newcommand{\Su}{S_{\mathbf u}}
\newcommand{\Lu}{L_{\mathbf u,\varphi}^{(N)}}
\newcommand{\Tg}{T_{\mathbf g}^{(n,N)}}
\newcommand{\dd}{\,d}
\newcommand{\ee}{\mathrm e}

\DeclareMathOperator{\dist}{dist}
\DeclareMathOperator{\Id}{Id}

\DeclareMathOperator{\Span}{span}

\begin{document}
\thispagestyle{firstpage}
\vspace*{20mm}
\begin{flushleft}
{\Large\bfseries
Carleson Measures and Integral\\[3pt]
Operators on Fock--Sobolev Spaces\par}
\end{flushleft}
\vspace{7mm}
\noindent{\large Jianxiang Dong$^{1}$ and Shibo Hu$^{1,*}$}
\par\vspace{3mm}
\begin{list}{}{\leftmargin=8mm\rightmargin=4mm}
\item[]\small\textbf{Abstract.}
For integers $n\ge1$ and $N\ge0$, let $\mathbf g=(g_0,\ldots,g_N)$
be a vector of entire functions. We consider the generalized integral operator
\[
 T_{\mathbf g}^{(n,N)}f
   =I^n\!\left(\sum_{j=0}^N g_jf^{(j)}\right),
 \qquad (Ih)(z)=\int_0^z h(w)\dd w.
\]
We characterize its boundedness and compactness between Fock--Sobolev
spaces $\F_m^p$ and $\F_\ell^q$, where $0<p,q<\infty$ and
$m,\ell\ge0$ are integers. The coefficients of critical total degree
determine the operator modulo compact operators and yield essential-norm
estimates. In the Hilbert case, a Laurent Toeplitz representation gives
the exact essential norm, a Fredholm criterion and an index formula.
As an application, we obtain conditions for solutions of linear
differential equations to belong to Fock--Sobolev spaces, including
a necessary and sufficient polynomial zero condition in the Hilbert case.
For entire functions $u_0,\ldots,u_N$ and $\varphi$, we also characterize
boundedness, compactness and essential norms of the operators
\[
 L_{\mathbf u,\varphi}^{(N)}f
   =\sum_{j=0}^N u_j\,f^{(j)}\circ\varphi.
\]
The proofs use Carleson measure characterizations for finite sums of
weighted differentiation operators on Fock--Sobolev spaces.
\end{list}

\section{Introduction and Main Results}
A Carleson criterion for a differential sum must account for cancellation between its terms. We use the coefficient polynomial of the sum to obtain a common measure criterion for integral and composition--differentiation operators on Fock--Sobolev spaces. The integral-operator criteria also lead to essential-norm, Fredholm and differential-equation results.

\subsection{Carleson Measures}
The scalar Carleson measure theory for Fock--Sobolev spaces was developed in \cite{CZ,MengCarleson}; for the classical Fock setting, see \cite{Zhu}. In \cite{AHTYY}, Arroussi, He, Tong, Yang and Yang introduced Carleson measures associated with finite sums of weighted differentiation operators on classical Fock spaces. We extend their characterization to Fock--Sobolev source spaces and obtain essential-norm estimates for all $0<p,q<\infty$.

Let $dA$ denote Lebesgue area measure on $\C$, and write $\Hc$ for the space of entire functions. Given $0<p<\infty$ and $m\in\Nzero$, define
\begin{equation}\label{eq:norm}
 \|f\|_{p,m}^{p}
 =\frac p{2\pi}\int_{\C}|f(z)|^p(1+|z|)^{mp}
                  \ee^{-p|z|^2/2}\dd A(z).
\end{equation}
We write $\F_m^p$ for the entire functions with finite norm or quasi-norm in \eqref{eq:norm}. This is the Fock--Sobolev space of order $m$ \cite{CZ}; in particular, $\F_0^p=\F^p$. We use $m,\ell\in\Nzero$ for the source and target Sobolev orders and $0<p,q<\infty$ for the corresponding exponents. Unless specified otherwise, $L^t$ means $L^t(\C,dA)$.

Fix $N\in\Nzero$, a vector $\mathbf u=(u_0,\ldots,u_N)$ with each $u_j\in\Hc$, and a positive locally finite Borel measure $\mu$ on $\C$. Define
\[
 \du f(z)=\sum_{j=0}^N u_j(z)f^{(j)}(z),\qquad
 P_{\mathbf u(z)}(\xi)=\sum_{j=0}^N u_j(z)\xi^j.
\]
In $P_{\mathbf u(z)}^{(k)}$, differentiation is with respect to $\xi$ while $z$ is fixed. We shall use the quantities
\begin{equation}\label{eq:jet}
 \begin{aligned}
 \Su(z)
   &=\sum_{k=0}^N|P_{\mathbf u(z)}^{(k)}(\bar z)|,\\[2pt]
 d\nu(z)
   &=\frac{\displaystyle\sum_{k=0}^N
       |P_{\mathbf u(z)}^{(k)}(\bar z)|^q}
       {(1+|z|)^{mq}}\dd\mu(z).
 \end{aligned}
\end{equation}
The measure $\nu$ is locally finite. For another finite coefficient vector $\mathbf v$, we use $P_{\mathbf v(z)}$ and $S_{\mathbf v}$ in the same way.
We call $\mu$ a $(\mathbf u,p,q;m)$-Carleson measure if there exists $C>0$ such that
\begin{equation}\label{eq:embedding}
 \int_{\C}|\du f(z)|^q\ee^{-q|z|^2/2}\dd\mu(z)
       \le C\|f\|_{p,m}^q,\qquad f\in\F_m^p.
\end{equation}
Denote this map into
$L^q(\C,\ee^{-q|z|^2/2}d\mu)$ by $J_{\mathbf u,\mu}f=\du f$. The notation $\|J_{\mathbf u,\mu}f\|_{L^q}$ always refers to this Gaussian-weighted measure. Such a Carleson measure is called vanishing if the left side of \eqref{eq:embedding} tends to zero along every bounded sequence in $\F_m^p$ that converges to zero uniformly on compact subsets of $\C$. By Lemma~\ref{lem:compact}, this is equivalent to compactness of $J_{\mathbf u,\mu}$.

For $\rho>0$, write $D(w,\rho)=\{z:|z-w|<\rho\}$ and
\[
 B_\rho\nu(w)=\nu(D(w,\rho)).
\]
Thus $B_\rho\nu$ is the disk-mass function, without normalization by the disk area. $(p,q)$-Fock Carleson measures and vanishing $(p,q)$-Fock Carleson measures correspond to the scalar case $m=N=0$ and $u_0=1$. We denote the essential norm of an operator $T$ by $\|T\|_{\mathrm e}$, with its source and target fixed.

\begin{theorem}\label{thm:carlesonup}
Assume $0<p\le q<\infty$, and let $\rho>0$. Then the following are equivalent:
\begin{enumerate}[ref=\alph*]
\item\label{cup:mu} The measure $\mu$ is $(\mathbf u,p,q;m)$-Carleson;
\item\label{cup:nu} The measure $\nu$ of \eqref{eq:jet} is a $(p,q)$-Fock Carleson measure;
\item\label{cup:mass} $\displaystyle\sup_{w\in\C}\nu(D(w,\rho))<\infty$.
\end{enumerate}
Whenever these conditions hold, the operator norm and essential norm satisfy
\begin{equation}\label{eq:essentialmeasure}
 \begin{aligned}
 \|J_{\mathbf u,\mu}\|^q
   &\asymp\|B_\rho\nu\|_\infty,
\qquad \|J_{\mathbf u,\mu}\|_{\mathrm e}^{q}
   &\asymp\limsup_{|w|\to\infty}B_\rho\nu(w).
 \end{aligned}
\end{equation}
Moreover, each of the following conditions is equivalent to compactness of $J_{\mathbf u,\mu}$:
\begin{enumerate}[start=4,ref=\alph*]
\item\label{cup:vanmu} $\mu$ is a vanishing $(\mathbf u,p,q;m)$-Carleson measure;
\item\label{cup:vannu} $\nu$ is a vanishing $(p,q)$-Fock Carleson measure;
\item\label{cup:vanmass} $\displaystyle\lim_{|w|\to\infty}\nu(D(w,\rho))=0$.
\end{enumerate}
\end{theorem}
We now turn to the complementary case \(0<q<p<\infty\). In this setting, the Carleson and vanishing Carleson conditions are equivalent, and both admit an equivalent characterization in terms of the \(L^r\)-integrability of \(B_\rho\nu\), as shown in the following theorem.
\begin{theorem}\label{thm:carlesondown}
Suppose $0<q<p<\infty$, and put $r=\frac{p}{p-q}$. For any $\rho>0$, the following conditions are equivalent:
\begin{enumerate}[ref=\alph*]
\item\label{cdown:mu} $\mu$ is a $(\mathbf u,p,q;m)$-Carleson measure;
\item\label{cdown:vanmu} $\mu$ is a vanishing $(\mathbf u,p,q;m)$-Carleson measure;
\item\label{cdown:nu} The measure $\nu$ of \eqref{eq:jet} is a $(p,q)$-Fock Carleson measure;
\item\label{cdown:vannu} $\nu$ is a vanishing $(p,q)$-Fock Carleson measure;
\item\label{cdown:mass} $B_\rho\nu\in L^r(\C,dA)$.
\end{enumerate}
Each of these conditions also characterizes compactness of $J_{\mathbf u,\mu}$. In this case,
\[
 \|J_{\mathbf u,\mu}\|^q\asymp\|B_\rho\nu\|_{L^r}.
\]
\end{theorem}
In either theorem, a disk-mass condition that holds for one radius $\rho>0$ holds for every positive radius, although the comparison constants may change. The measure in \eqref{eq:jet} involves $P_{\mathbf u(z)}^{(k)}(\bar z)$ for every $0\le k\le N$. Setting $m=0$ or $N=0$ recovers, respectively, the differential criterion of \cite{AHTYY} or the scalar Fock--Sobolev criterion.

\subsection{Generalized Integral Operators}
Volterra operators on Fock spaces were studied in \cite{Constantin}, with spectral results on Fock--Sobolev spaces in \cite{MengSpectrum} and generalized operators in \cite{AbdellaWorku}. The finite sums in \cite{AHTYY} motivate the following extension, in which the integration order, differential order and two Sobolev orders may vary independently.

Let $D$ denote differentiation. For $n\in\N$, $N\in\Nzero$ and an entire coefficient vector $\mathbf g=(g_0,\ldots,g_N)$, set
\[
 (Ih)(z)=\int_0^z h(w)\dd w,\qquad
 \Tg f=I^n\!\left(\sum_{j=0}^N g_jf^{(j)}\right).
\]
The first $n$ Taylor coefficients of $I^n h$ vanish. The differential order $N$ is independent of $n$; the case $N=n-1$ includes the sums in \cite{AHTYY}, while $N\ge n$ also allows companion terms. Fix source and target orders $m,\ell$, and put $d=n+m-\ell$. We set $\deg0=-\infty$ and assign degree $+\infty$ to a nonpolynomial entire function.

\begin{theorem}\label{thm:integral}
Let $T=\Tg$ act from $\F_m^p$ to $\F_\ell^q$.
\begin{enumerate}
\item Suppose $p\le q$. Then $T$ is bounded precisely when $\deg g_j+j\le d$ for every $0\le j\le N$, and compact precisely when all these inequalities are strict.
\item Suppose $q<p$. Then boundedness and compactness are equivalent, and they hold precisely when
\begin{equation}\label{eq:degree-down}
 \deg g_j+j<d-2\big(\frac1q-\frac1p\big),\qquad 0\le j\le N.
\end{equation}
\end{enumerate}
\end{theorem}
For $m=\ell=0$ and $N=n-1$, this recovers the degree conditions in \cite{AHTYY}. Equal Sobolev orders leave those thresholds unchanged. The inequality in \eqref{eq:degree-down} remains strict when the threshold is an integer; the number $2$ comes from planar area measure.

The operator introduced by Chalmoukis \cite{Chalmoukis} is a particular instance:
\[
 T_{g,\mathbf a}f=I^n\!\left(\sum_{j=0}^{n-1}
                   a_jg^{(n-j)}f^{(j)}\right),
 \qquad \mathbf a\in\C^n.
\]
The degree threshold depends on the largest index with a nonzero coefficient.

\begin{corollary}\label{cor:chalmoukis}
Assume $\mathbf a\ne0$, and set
\[
 j_* =\max\{j:a_j\ne0\},\qquad h=n-j_*,\qquad \Delta=m-\ell.
\]
Then $T_{g,\mathbf a}=0$ exactly when $\deg g<h$. If $\deg g\ge h$, the criteria are as follows:
\begin{enumerate}
\item For $p\le q$, $T_{g,\mathbf a}$ is bounded exactly when $\deg g\le2h+\Delta$, and compact exactly when $\deg g<2h+\Delta$;
\item For $q<p$, it is bounded (equivalently, compact) exactly when
$\deg g<2h+\Delta-2\big(\frac1q-\frac1p\big)$.
\end{enumerate}
Consequently, $g$ must be a polynomial whenever $T_{g,\mathbf a}$ is bounded and nonzero.
If $\mathbf a=0$, then $T_{g,\mathbf a}=0$ for every entire $g$.
\end{corollary}

In the case $a_{n-1}\ne0$ and $m=\ell$, boundedness for $p\le q$ allows symbols of degree at most two, while compactness requires degree at most one. The assumption $a_{n-1}\ne0$ is essential to this particular statement.

Essential norms of individual integral operators have been studied in \cite{MengEssential,WorkuYesuf}. For the finite sums above, the coefficients at the critical total degree determine the operator modulo compact operators. We encode them in a Laurent polynomial.

Let $\mathbb T=\{\zeta:|\zeta|=1\}$ and $\mathbb D=\{w:|w|<1\}$. For a bounded $T$ with $p\le q$, define
\[
 \mathcal I_d=\{j\in\Nzero:0\le j\le N,\ j\le d\},
 \qquad c_j=[z^{d-j}]g_j\quad(j\in\mathcal I_d).
\]
The symbol $[z^h]g$ denotes the coefficient of $z^h$ in $g$. The critical Laurent polynomial is
\begin{equation}\label{eq:criticalsymbol}
 \sigma_{\mathbf g}(\zeta)
       =\sum_{j\in\mathcal I_d}c_j\zeta^{n+d-2j}.
\end{equation}
We understand an empty sum to be zero; exponents in \eqref{eq:criticalsymbol} may be negative.

In the Hilbert setting, let $(e_k)_{k\ge0}$ be the canonical basis of $\ell^2(\Nzero)$ and define the unitary map $U_m:\ell^2(\Nzero)\to\F_m^2$ by $U_m e_k=z^k/\|z^k\|_{2,m}$. With $S e_k=e_{k+1}$, write $V_s=S^s$ for $s\ge0$ and $V_s=(S^*)^{-s}$ for $s<0$.

\begin{theorem}\label{thm:essentialintegral}
Let $0<p\le q<\infty$, and suppose that $T=\Tg:\F_m^p\to\F_\ell^q$ is bounded. Then
\begin{equation}\label{eq:essentialintegral}
                  \|T\|_{\mathrm e}\asymp\max_{\mathbb T}|\sigma_{\mathbf g}|.
\end{equation}
In addition, $T$ differs by a compact operator from
\[
 f\longmapsto I^n\!\left(\sum_{j\in\mathcal I_d}c_jz^{d-j}f^{(j)}\right).
\]
With $n,N,m,\ell,p,q$ fixed, the bounded operators of this form, modulo compact operators of the same form, constitute a vector space of dimension $|\mathcal I_d|$. Two bounded operators of this form have compact difference if and only if their critical coefficients coincide.

For $p=q=2$, the operator has the representation
\begin{equation}\label{eq:toeplitzmodel}
 \begin{aligned}
 U_\ell^{-1}TU_m
   &=A_{\sigma_{\mathbf g}}+K,
 \qquad A_{\sigma_{\mathbf g}}
   &=\sum_{j\in\mathcal I_d}c_jV_{n+d-2j},
 \end{aligned}
\end{equation}
with $K$ compact. Its essential norm is exactly
\begin{equation}\label{eq:essentialexact}
         \|T:\F_m^2\to\F_\ell^2\|_{\mathrm e}
                         =\max_{\mathbb T}|\sigma_{\mathbf g}|.
\end{equation}
\end{theorem}
The identity \eqref{eq:essentialexact} uses the normalization in \eqref{eq:norm}; under equivalent norms one retains the comparison in \eqref{eq:essentialintegral}. When $q<p$, Theorem~\ref{thm:integral} makes every bounded $T$ compact, so its essential norm is zero.

Fulsche \cite{Fulsche} connected rotationally continuous operators on classical Fock space with band-dominated matrices in the monomial basis. In the present setting, \eqref{eq:toeplitzmodel} describes the limiting diagonals even when the Sobolev orders differ. Its Laurent symbol also governs Fredholmness.

For a continuous function $a$ without zeros on the positively oriented circle $\mathbb T$, let $\operatorname{wind}(a,0)$ be its winding number about the origin. We take $\operatorname{ind}T=\dim\ker T-\operatorname{codim}\operatorname{ran}T$. For an operator on a single space, its Fredholm essential spectrum is $\sigma_{\mathrm e}(T)=\{\lambda:T-\lambda\Id\text{ is not Fredholm}\}$.

\begin{theorem}\label{thm:fredholmsymbol}
Suppose $T=\Tg:\F_m^2\to\F_\ell^2$ is bounded. Then
\begin{equation}\label{eq:fredholmcriterion}
 T\text{ is Fredholm}
       \quad\Longleftrightarrow\quad
 \sigma_{\mathbf g}(\zeta)\ne0\quad(\zeta\in\mathbb T).
\end{equation}
In that case,
\begin{equation}\label{eq:indexsymbol}
             \operatorname{ind}T=-\operatorname{wind}(\sigma_{\mathbf g},0).
\end{equation}
If $m=\ell$, we also have
\begin{equation}\label{eq:essspectrum}
 \begin{aligned}
 \sigma_{\mathrm e}(T)
   &=\sigma_{\mathbf g}(\mathbb T),\\
 \operatorname{ind}(T-\lambda\Id)
   &=-\operatorname{wind}(\sigma_{\mathbf g}-\lambda,0).
 \end{aligned}
\end{equation}
The index formula in \eqref{eq:essspectrum} holds for every $\lambda\notin\sigma_{\mathbf g}(\mathbb T)$.
\end{theorem}
The Fredholm theorem for continuous Toeplitz symbols is classical; see \cite{Douglas}. After establishing \eqref{eq:toeplitzmodel}, we give a proof for the finite Laurent polynomial at hand to track negative shifts and the sign of the index. The essential-spectrum statement in \eqref{eq:essspectrum} concerns an operator on one Fock--Sobolev space and makes no identification between distinct spaces.

\subsection{Linear Differential Equations}
Consider the linear differential equation
\begin{equation}\label{eq:ode}
 f^{(n)}+\sum_{j=0}^{n-1}g_jf^{(j)}=G,
\end{equation}
where the coefficients and the forcing term are entire. Related conditions in weighted Fock spaces were obtained in \cite{HuODE}. Assume that
\begin{equation}\label{eq:odecoeff}
 g_j(z)=c_jz^{n-j}+r_j(z),\qquad
 \deg r_j<n-j,\quad 0\le j<n.
\end{equation}
Write $T=T_{\mathbf g}^{(n,n-1)}$. By Theorem~\ref{thm:integral}, this operator is bounded on $\F_m^p$ for every $m\in\Nzero$ and $0<p<\infty$.

\begin{theorem}\label{thm:ode}
Assume $0<p<\infty$, $m\in\Nzero$ and $I^nG\in\F_m^p$. Under \eqref{eq:odecoeff}, put $\theta=\min\{1,p\}$. There is a constant $\varepsilon_{n,m,p}>0$ such that every entire solution of \eqref{eq:ode} lies in $\F_m^p$ provided that
\[
 \left(\sum_{j=0}^{n-1}|c_j|^\theta\right)^{1/\theta}
                          <\varepsilon_{n,m,p}.
\]
No restriction is imposed on the lower degree coefficients in $r_j$. In particular, the assertion applies when every $c_j$ vanishes.
\end{theorem}

In the Hilbert case, the criterion can be stated exactly. Define
\begin{equation}\label{eq:odepolynomial}
                   Q(w)=1+\sum_{j=0}^{n-1}c_jw^{n-j}.
\end{equation}
For a polynomial $Q$, denote by $Z_{\mathbb D}(Q)$ the number of its zeros in $\mathbb D$, counted with multiplicity.

\begin{theorem}\label{thm:odehilbert}
Assume \eqref{eq:odecoeff}, let $m\in\Nzero$, and define $Q$ by \eqref{eq:odepolynomial}. The operator $\Id+T$ is Fredholm on $\F_m^2$ exactly when $Q$ has no zeros on $\mathbb T$. In this case,
\begin{equation}\label{eq:odeindex}
              \operatorname{ind}(\Id+T)=-2Z_{\mathbb D}(Q).
\end{equation}
The following three statements are equivalent:
\begin{enumerate}
\item $\Id+T$ has a bounded inverse on $\F_m^2$;
\item The polynomial $Q$ is zero-free on $\overline{\mathbb D}$;
\item For each entire $G$ satisfying $I^nG\in\F_m^2$ and each prescribed set of values $f^{(j)}(0)$, $0\le j<n$, the entire solution of \eqref{eq:ode} lies in $\F_m^2$.
\end{enumerate}
When these conditions hold, every such solution satisfies
\[
 \|f\|_{2,m}\le C\left(\|I^nG\|_{2,m}
                     +\sum_{j=0}^{n-1}|f^{(j)}(0)|\right).
\]
Here $C$ may depend on $n,m$ and all the coefficients $g_j$. The zero condition does not depend on $m$ or on the lower degree terms $r_j$.
\end{theorem}
\subsection{Sums of Weighted Composition--Differentiation Operators}
Weighted composition operators on Fock and Fock--Sobolev spaces have been studied in \cite{TienKhoi,HeCao}. As in \cite{AHTYY}, we consider finite sums whose terms share one composition symbol. Let $C_\varphi f=f\circ\varphi$. For $\varphi\in\Hc$, put
\[
 \Lu f(z)=\sum_{j=0}^N u_j(z)f^{(j)}(\varphi(z)).
\]
We shall use the function
\begin{equation}\label{eq:B}
 \mathcal B(z)=\frac{(1+|z|)^\ell}{(1+|\varphi(z)|)^m}
 \ee^{(|\varphi(z)|^2-|z|^2)/2}
 \sum_{k=0}^N|P_{\mathbf u(z)}^{(k)}(\overline{\varphi(z)})|.
\end{equation}
The weight factors in the numerator and denominator of \eqref{eq:B} arise from the target and source Fock--Sobolev norms, respectively.

\begin{proposition}\label{prop:affine}
Assume $\mathbf u\not\equiv0$. Boundedness of $\Lu:\F_m^p\to\F_\ell^q$ forces $\varphi(z)=az+b$ with $|a|\le1$. If $a=0$, the operator is bounded exactly when it is compact, and this occurs exactly when $u_j\in\F_\ell^q$ for every $j$. No restriction on $\varphi$ is needed when $\mathbf u\equiv0$.
\end{proposition}

\begin{theorem}\label{thm:compup}
Assume $0<p\le q<\infty$ and $\varphi(z)=az+b$ with $0<|a|\le1$. Then:
\begin{enumerate}
\item Boundedness of $\Lu:\F_m^p\to\F_\ell^q$ is equivalent to
\[
                  \sup_{z\in\C}\mathcal B(z)<\infty.
\]
\item Compactness of $\Lu:\F_m^p\to\F_\ell^q$ is equivalent to $\mathcal B(z)\to0$ as $|z|\to\infty$.
\end{enumerate}
Whenever $\Lu$ is bounded, its operator norm and essential norm satisfy
\[
 \|\Lu\|\asymp\|\mathcal B\|_{L^\infty},\qquad
 \|\Lu\|_{\mathrm e}\asymp\limsup_{|z|\to\infty}\mathcal B(z).
\]
The constants in these comparisons may depend on $a$ as well as on $p,q,m,\ell,N$.
\end{theorem}

\begin{theorem}\label{thm:compdown}
Assume $0<q<p<\infty$, put $s=\frac{pq}{p-q}$, and let $\varphi(z)=az+b$ with $0<|a|\le1$. Then the following are equivalent:
\begin{enumerate}
\item $\Lu:\F_m^p\to\F_\ell^q$ is bounded;
\item $\Lu:\F_m^p\to\F_\ell^q$ is compact;
\item $\mathcal B\in L^s(\C,dA)$.
\end{enumerate}
Under these conditions, $\|\Lu\|\asymp\|\mathcal B\|_{L^s}$ and $\|\Lu\|_{\mathrm e}=0$. The comparison constants may depend on the fixed affine map.
\end{theorem}

\begin{corollary}\label{cor:commute}
For each $n\in\N$, $D^nC_\varphi:\F_m^p\to\F_\ell^q$ is bounded (respectively, compact) if and only if $C_\varphi D^n:\F_m^p\to\F_\ell^q$ is bounded (respectively, compact).
\end{corollary}

Section~\ref{sec:preliminaries} collects derivative estimates, test functions and compact approximation results. In Section~\ref{sec:carleson} we prove Theorems~\ref{thm:carlesonup} and \ref{thm:carlesondown}. Section~\ref{sec:operators} then treats generalized integral operators, essential norms and Fredholm properties, differential equations, and sums of weighted composition--differentiation operators.

The constants implicit in $\lesssim$ and $\asymp$ may depend on fixed orders, exponents and radii, but are independent of the coefficient vectors and measures. We indicate any dependence on a fixed affine map.

\section{Preliminaries}\label{sec:preliminaries}
We collect the estimates needed for the Carleson criteria and their operator consequences. Throughout, the integrability exponents $p$ and $q$ are positive and finite.

\subsection{Derivative Estimates and Test Functions}
By the derivative characterization of Cho and Zhu \cite{CZ},
\begin{equation}\label{eq:sobolevnorm}
 \|f\|_{p,m}\asymp
 \sum_{j=0}^{m-1}|f^{(j)}(0)|+\|f^{(m)}\|_{p,0}.
\end{equation}
The sum is empty when $m=0$. Finiteness of either side is equivalent to membership in $\F_m^p$.

We also use the classical higher-derivative formula \cite{HuNorm,Ueki}:
\begin{equation}\label{eq:classicalLP}
 \|h\|_{q,0}^q\asymp\sum_{j=0}^{v-1}|h^{(j)}(0)|^q
   +\int_{\C}\frac{|h^{(v)}(z)|^q}{(1+|z|)^{vq}}
                                  \ee^{-q|z|^2/2}\dd A(z),
\end{equation}
for $v\in\N$, $0<q<\infty$ and $h\in\Hc$, with membership in $\F^q$ equivalent to finiteness of the right-hand side. The finite Taylor sum is needed also when $q<1$.

\begin{lemma}\label{lem:LP}
Let $n\in\N$, $\ell\in\Nzero$, and $H\in\Hc$. Then
\begin{equation}\label{eq:primitive}
 \|I^nH\|_{q,\ell}^q\asymp
 \int_{\C}|H(z)|^q(1+|z|)^{(\ell-n)q}
                                    \ee^{-q|z|^2/2}\dd A(z).
\end{equation}
The equivalence includes the case in which both sides are infinite.
\end{lemma}
\begin{proof}
Since $(I^nH)^{(j)}(0)=0$ for $0\le j<n$, the case $\ell<n$ follows by applying \eqref{eq:sobolevnorm} to $I^nH$ and \eqref{eq:classicalLP} to $I^{n-\ell}H$, with $v=n-\ell$:
\[
 \|I^nH\|_{q,\ell}^q\asymp\|I^{n-\ell}H\|_{q,0}^q
 \asymp\int_{\C}|H(z)|^q(1+|z|)^{(\ell-n)q}
                  \ee^{-q|z|^2/2}\dd A(z).
\]
For $\ell=n$, \eqref{eq:sobolevnorm} reduces the left-hand side to $\|H\|_{q,0}^q$. If $\ell>n$, set $v=\ell-n$. The vanishing derivatives below order $n$, together with $(I^nH)^{(n+j)}(0)=H^{(j)}(0)$ and $(I^nH)^{(\ell)}=H^{(v)}$, allow \eqref{eq:sobolevnorm} to be applied first to $I^nH$ and then to $H$:
\[
 \|I^nH\|_{q,\ell}^q
 \asymp\sum_{j=0}^{v-1}|H^{(j)}(0)|^q+\|H^{(v)}\|_{q,0}^q
 \asymp\|H\|_{q,v}^q.
\]
By the definition of $\|H\|_{q,v}$, the last expression is comparable to the integral on the right of \eqref{eq:primitive}. For a fixed finite number of nonnegative terms, the $q$th power of their sum is comparable to the sum of their $q$th powers for every $q>0$. The cited norm characterizations also give equivalence of finiteness, so the comparison includes infinite values.
\end{proof}

Given $f\in\Hc$ and $k\in\Nzero$, write
\begin{equation}\label{eq:Q}
 Q_kf(z)=\ee^{-|z|^2/2}
     \sum_{j=0}^k\binom{k}{j}(-\bar z)^{k-j}f^{(j)}(z).
\end{equation}
In this notation $\bar z$ is treated as a fixed parameter. The local estimate below is obtained from an entire function of the translation variable.

\begin{lemma}\label{lem:local}
If $0<p<\infty$, $m,k\in\Nzero$, and $\rho>0$, then
\begin{equation}\label{eq:local}
 (1+|z|)^{mp}|Q_kf(z)|^p
 \lesssim\int_{D(z,\rho)}|f(\zeta)|^p(1+|\zeta|)^{mp}
                       \ee^{-p|\zeta|^2/2}\dd A(\zeta).
\end{equation}
Consequently,
\begin{equation}\label{eq:Qglobal}
 \|(1+|\cdot|)^m Q_kf\|_{L^p}\lesssim\|f\|_{p,m},
\end{equation}
and
\begin{equation}\label{eq:point}
 |f^{(k)}(z)|\ee^{-|z|^2/2}
      \lesssim(1+|z|)^{k-m}\|f\|_{p,m}.
\end{equation}
Also $\F_m^p$ embeds continuously in $\F_m^q$ whenever $p\le q$.
\end{lemma}
\begin{proof}
For fixed $z$, the function
\[
 H_z(\eta)=f(z+\eta)\ee^{-\bar z\eta-|z|^2/2}
\]
is entire, and
\[
 H_z^{(k)}(0)=Q_kf(z),\qquad
 |H_z(\eta)|\ee^{-|\eta|^2/2}
      =|f(z+\eta)|\ee^{-|z+\eta|^2/2}.
\]
Cauchy's estimate on $D(0,\rho/2)$, followed by the submean inequality for $|H_z|^p$, gives
\[
 |H_z^{(k)}(0)|^p
 \lesssim\sup_{|\eta|\le\rho/2}|H_z(\eta)|^p
 \lesssim\int_{D(0,\rho)}|H_z(\eta)|^p\dd A(\eta).
\]
On $D(0,\rho)$ the factor $\ee^{p|\eta|^2/2}$ is bounded, while $1+|z+\eta|\asymp1+|z|$. Changing variables to $\zeta=z+\eta$ proves \eqref{eq:local}. Integrating that estimate over $z$ and applying Tonelli's theorem gives \eqref{eq:Qglobal}. Binomial inversion in \eqref{eq:Q} yields
\begin{equation}\label{eq:inverseQ}
 f^{(k)}(z)\ee^{-|z|^2/2}
       =\sum_{j=0}^k\binom{k}{j}\bar z^{k-j}Q_jf(z).
\end{equation}
Indeed, after substituting \eqref{eq:Q}, the coefficient of $f^{(v)}(z)\ee^{-|z|^2/2}$ on the right is
\[
 \sum_{j=v}^k\binom{k}{j}\binom{j}{v}(-1)^{j-v}\bar z^{k-v}
 =\binom{k}{v}\bar z^{k-v}(1-1)^{k-v}.
\]
This coefficient is $1$ when $v=k$ and $0$ otherwise. Since \eqref{eq:local} gives $|Q_jf(z)|\lesssim(1+|z|)^{-m}\|f\|_{p,m}$, substitution into \eqref{eq:inverseQ} proves \eqref{eq:point}.

Set $F(z)=|f(z)|(1+|z|)^m\ee^{-|z|^2/2}$. The case $k=0$ of \eqref{eq:local} shows that $\|F\|_\infty\lesssim\|f\|_{p,m}$. Hence, when $p\le q$,
\[
 \int_{\C} F^q\dd A\le\|F\|_\infty^{q-p}
                 \int_{\C}F^p\dd A\lesssim\|f\|_{p,m}^q.
\]
This proves the continuous inclusion $\F_m^p\hookrightarrow\F_m^q$. The submean argument applies for all $p>0$.
\end{proof}

For $w\in\C$ and $k\in\{0,\ldots,N\}$, introduce
\begin{equation}\label{eq:atoms}
 k_w(z)=\ee^{\bar wz-|w|^2/2},\qquad
 E_{w,k}(z)=\frac{(z-w)^k k_w(z)}{(1+|w|)^m}.
\end{equation}
The functions $E_{w,k}$ are test functions, not normalized reproducing kernels of $\F_m^2$.

\begin{lemma}\label{lem:atoms}
Uniformly in their centers, the test functions obey
\[
 \sup_{w\in\C}\|E_{w,k}\|_{p,m}<\infty.
\]
They converge to zero uniformly on compact sets as $|w|\to\infty$.
Moreover, for $\Lambda=\rho(\mathbb Z+i\mathbb Z)$ and every finitely supported $c=(c_w)_{w\in\Lambda}$,
\begin{equation}\label{eq:synthesis}
 \left\|\sum_{w\in\Lambda}c_wE_{w,k}\right\|_{p,m}
                  \lesssim\|c\|_{\ell^p}.
\end{equation}
The implicit constant is independent of $c$ and of its support.
\end{lemma}
\begin{proof}
The identity $|k_w(z)|\ee^{-|z|^2/2}=\ee^{-|z-w|^2/2}$, together with the inequality $1+|z|\le(1+|w|)(1+|z-w|)$, gives
\begin{align*}
 |E_{w,k}(z)|(1+|z|)^m\ee^{-|z|^2/2}
 &=|z-w|^k\frac{(1+|z|)^m}{(1+|w|)^m}
                         \ee^{-|z-w|^2/2}\\
 &\le(1+|z-w|)^{m+k}\ee^{-|z-w|^2/2}
 =:H(z-w).
\end{align*}
Because $H\in L^p(\C,dA)$ for every $p>0$, the first assertion about the norms follows. On $|z|\le R$, the unweighted modulus of $E_{w,k}$ is at most
\[
 (R+|w|)^k(1+|w|)^{-m}
                    \ee^{R|w|-|w|^2/2}.
\]
The Gaussian decay in $|w|$ makes this bound tend to zero uniformly for $|z|\le R$.

The lattice separation also gives $\sup_z\sum_{w\in\Lambda}H(z-w)<\infty$: in the annulus $j\le|z-w|<j+1$ there are at most $C_\rho(1+j)^2$ lattice points, and the corresponding Gaussian series converges uniformly in $z$. For $0<p\le1$, subadditivity of the $p$th power yields
\[
 \int_{\C}\left(\sum_w|c_w|H(z-w)\right)^p\dd A(z)
 \le\sum_w|c_w|^p\int_{\C}H(z-w)^p\dd A(z)
 \lesssim\sum_w|c_w|^p.
\]
For $p>1$, discrete H\"older's inequality, with weights $H(z-w)$, instead gives
\begin{align*}
 \left(\sum_w|c_w|H(z-w)\right)^p
 &\le\left(\sum_w|c_w|^pH(z-w)\right)
       \left(\sum_wH(z-w)\right)^{p-1}\\
 &\lesssim\sum_w|c_w|^pH(z-w).
\end{align*}
Integrating in $z$ and using translation invariance, we obtain
\[
 \int_{\C}\left(\sum_w|c_w|H(z-w)\right)^p\dd A(z)
       \lesssim\|H\|_{L^1}\sum_w|c_w|^p.
\]
Combining the two ranges establishes \eqref{eq:synthesis}.
\end{proof}

\subsection{Coefficient Polynomials}
The coefficient polynomial records cancellation between the terms of the differential operator. Its derivatives satisfy the local test estimate needed for the measure in \eqref{eq:jet}.

\begin{lemma}\label{lem:algebra}
Every $f\in\Hc$ satisfies
\begin{equation}\label{eq:algebra}
 \du f(z)\ee^{-|z|^2/2}
       =\sum_{k=0}^N\frac{P_{\mathbf u(z)}^{(k)}(\bar z)}{k!}
                       Q_kf(z).
\end{equation}
Given $\rho>0$, the following comparison is uniform for $|z-w|<\rho$:
\begin{equation}\label{eq:localtests}
 \sum_{k=0}^N|\du E_{w,k}(z)|^q\ee^{-q|z|^2/2}
 \asymp\frac{\displaystyle\sum_{k=0}^N
                   |P_{\mathbf u(z)}^{(k)}(\bar z)|^q}
                  {(1+|z|)^{mq}}.
\end{equation}
\end{lemma}
\begin{proof}
Substituting \eqref{eq:inverseQ} into the differential sum and collecting the coefficients of $Q_kf$ gives
\begin{align*}
 \du f(z)\ee^{-|z|^2/2}
 &=\sum_{j=0}^N u_j(z)
          \sum_{k=0}^j\binom{j}{k}\bar z^{j-k}Q_kf(z)\\
 &=\sum_{k=0}^N Q_kf(z)
          \sum_{j=k}^N\binom{j}{k} u_j(z)\bar z^{j-k}\\
 &=\sum_{k=0}^N\frac{1}{k!}
                 P_{\mathbf u(z)}^{(k)}(\bar z)Q_kf(z).
\end{align*}
The interchange is legitimate because both sums are finite; the last line is \eqref{eq:algebra}.

To compare the two sides of \eqref{eq:localtests}, fix $z,w$ with $|z-w|<\rho$ and write $x=z-w$ and $v_j=P_{\mathbf u(z)}^{(j)}(\bar w)$. Leibniz's rule, applied for $0\le i,k\le N$, yields
\[
 \frac{d^i}{dz^i}\big((z-w)^k k_w(z)\big)
 =k_w(z)\sum_{j=0}^{\min\{i,k\}}
       \binom{i}{j}\frac{k!}{(k-j)!}x^{k-j}\bar w^{i-j}.
\]
Multiplication by $u_i(z)$ and summation in $i$ therefore yield
\begin{align*}
 \du E_{w,k}(z)
 &=\frac{k_w(z)}{(1+|w|)^m}
   \sum_{j=0}^k\frac{k!}{(k-j)!}x^{k-j}
      \sum_{i=j}^N\binom{i}{j} u_i(z)\bar w^{i-j}\\
 &=\frac{k_w(z)}{(1+|w|)^m}A_k,
\end{align*}
where
\[
 A_k=\sum_{j=0}^k\binom{k}{j} x^{k-j}v_j.
\]
The binomial inversion formula reads
\[
 v_k=\sum_{j=0}^k\binom{k}{j}(-x)^{k-j}A_j.
\]
Because $|x|<\rho$, the coefficients in both triangular systems are bounded in terms of $N$ and $\rho$. Apply $|\sum_{j=0}^N a_j|^q\le (N+1)^{\max\{q-1,0\}}\sum_{j=0}^N|a_j|^q$ to each system to obtain
\[
 C^{-1}\sum_{j=0}^N|v_j|^q
       \le\sum_{k=0}^N|A_k|^q
       \le C\sum_{j=0}^N|v_j|^q.
\]
\Needspace{8\baselineskip}
It remains to replace $\bar w$ by $\bar z$. Since the polynomial variable is independent of the coefficient variable, Taylor's formula gives
\[
 P_{\mathbf u(z)}^{(k)}(\bar z)
 =\sum_{j=k}^N\frac{v_j}{(j-k)!}(\bar z-\bar w)^{j-k},
\]
\nopagebreak
and
\nopagebreak
\[
 v_k=\sum_{j=k}^N
       \frac{P_{\mathbf u(z)}^{(j)}(\bar z)}{(j-k)!}
                                     (\bar w-\bar z)^{j-k}.
\]
As $|z-w|<\rho$, the same finite-sum estimate applied in both directions shows that $\sum_k|v_k|^q\asymp\sum_k|P_{\mathbf u(z)}^{(k)}(\bar z)|^q$. Moreover,
\[
 |k_w(z)|^q\ee^{-q|z|^2/2}=\ee^{-q|z-w|^2/2},
 \qquad 1+|z|\asymp1+|w|\quad(|z-w|<\rho).
\]
The Gaussian factor is uniformly bounded above and below on $|z-w|<\rho$. These comparisons prove \eqref{eq:localtests}; their constants do not involve $z,w$ or the coefficient functions.
\end{proof}

For the operator applications, we need a submean estimate for the derivatives of the coefficient polynomial evaluated along an affine map. This estimate allows us to pass from disk-mass bounds to pointwise estimates.

\begin{lemma}\label{lem:submeanjet}
Suppose that the entries of $\mathbf v=(v_0,\ldots,v_N)$ are entire and that $\eta(z)=az+b$. Let $W:\C\to(0,\infty)$ be measurable and locally comparable in the following sense: for some $C_W\ge1$ independent of $z,w$,
\[
 C_W^{-1}\le\frac{W(z)}{W(w)}\le C_W\quad\text{if }|z-w|<1.
\]
Suppose
\[
 \beta(z)=\frac\gamma2|z|^2+\operatorname{Re}(\lambda z)+\beta_0,
 \qquad \gamma,\beta_0\in\mathbb R,\quad\lambda\in\C.
\]
Put
\[
 A(z)=W(z)\ee^{\beta(z)}
           \sum_{k=0}^N|P_{\mathbf v(z)}^{(k)}(\overline{\eta(z)})|.
\]
Then, for any $t>0$ and fixed $\rho>0$,
\begin{equation}\label{eq:jetsubmean}
                 A(w)^t\lesssim\int_{D(w,\rho)}A(z)^t\dd A(z).
\end{equation}
The implicit constant may depend on $t,N,\rho,a,\gamma,C_W$, but is uniform in $w$ and $\mathbf v$. The inequality includes the case of an infinite integral.
\end{lemma}
\begin{proof}
The local comparability assumed for $W$ extends to every fixed radius: divide the segment from $w$ to $z$ into a bounded number of pieces of length less than one. Thus $W(z)\asymp W(w)$ whenever $|z-w|<\rho$. For fixed $w$, put $L_w=\gamma\bar w+\lambda$ and consider the entire function
\[
 H_k(z)=P_{\mathbf v(z)}^{(k)}(\overline{\eta(w)})
                         \ee^{L_w(z-w)}.
\]
Its modulus to the $t$th power is subharmonic, so the disk submean inequality applies. We also have
\[
 \beta(z)-\beta(w)=\operatorname{Re}(L_w(z-w))
                              +\frac\gamma2|z-w|^2,
\]
which, together with the submean inequality, gives
\begin{align*}
 &|P_{\mathbf v(w)}^{(k)}(\overline{\eta(w)})|^t
       \ee^{t\beta(w)}
 \lesssim\int_{D(w,\rho)}
       |P_{\mathbf v(z)}^{(k)}(\overline{\eta(w)})|^t
            \ee^{t\beta(z)}\dd A(z).
\end{align*}
The factor involving $\gamma|z-w|^2$ is bounded above and below on $D(w,\rho)$. Expand the polynomial about $\overline{\eta(z)}$ to obtain
\[
 P_{\mathbf v(z)}^{(k)}(\overline{\eta(w)})
 =\sum_{j=k}^N\frac{P_{\mathbf v(z)}^{(j)}(\overline{\eta(z)})}{(j-k)!}
                   \big(\overline{\eta(w)-\eta(z)}\big)^{j-k}.
\]
Since $|\eta(w)-\eta(z)|\le|a|\rho$ on the disk, the finite sum is bounded by a constant times the jet sum at $z$. Multiply by $W(w)^t$, use its local comparison with $W(z)^t$, and sum in $k$. The equivalence of $(\sum_{k=0}^N x_k)^t$ and $\sum_{k=0}^Nx_k^t$ for finitely many $x_k\ge0$ proves \eqref{eq:jetsubmean}, including for $t<1$.
\end{proof}

\subsection{Compact Approximation}
The following lemmas give a sequential compactness criterion and the compact approximants needed for the essential-norm estimates.

\begin{lemma}\label{lem:compact}
Let $T$ be either a bounded differential embedding as above, or a bounded linear operator into a Fock--Sobolev space whose formula preserves locally uniform convergence of entire functions. Then $T$ is compact if and only if every bounded sequence converging to zero uniformly on compact sets has images tending to zero in the target quasi-norm.
\end{lemma}
\begin{proof}
Local boundedness follows from Lemma~\ref{lem:local}. Montel's theorem and Fatou's lemma give a locally convergent subsequence of every bounded source sequence, with limit in $\F_m^p$. Subtracting this limit proves sufficiency.

Conversely, compactness gives a target-convergent subsequence of every image sequence. Point evaluations identify its limit with zero for the holomorphic target. For an $L^q$ target, a further subsequence converges almost everywhere, which gives the same conclusion by local convergence of derivatives. Every subsequence thus has a further subsequence whose image norm tends to zero, proving necessity.
\end{proof}
Completeness of $\F_m^p$ follows from the local estimate and completeness of weighted $L^p$.

We next construct compact dilations and establish the uniform tail estimates used in the essential-norm proofs.
\begin{lemma}\label{lem:dilations}
The dilation $R_tf(z)=f(tz)$ is compact on $\F_m^p$ for each $0<t<1$. Moreover,
\[
 \sup_{1/2\le t<1}\|R_t\|<\infty.
\]
For any fixed $R>0$ and $k\in\Nzero$, one also has
\begin{equation}\label{eq:dilationlocal}
 \sup_{\|f\|_{p,m}\le1}\sup_{|z|\le R}
       |((\Id-R_t)f)^{(k)}(z)|\longrightarrow0\quad(t\longrightarrow1).
\end{equation}
\end{lemma}
\begin{proof}
Under the change of variables $w=tz$, the elementary bound $1+|w|/t\le t^{-1}(1+|w|)$ gives
\begin{align*}
 \|R_tf\|_{p,m}^p
 &=\frac p{2\pi t^2}\int_{\C}|f(w)|^p(1+|w|/t)^{mp}
                     \ee^{-p|w|^2/(2t^2)}\dd A(w)\\
 &\le t^{-2-mp}\|f\|_{p,m}^p.
\end{align*}
Applying \eqref{eq:point} with $k=0$ at $tz$ and comparing the polynomial weights yields
\[
 |f(tz)|(1+|z|)^m\ee^{-|z|^2/2}
 \lesssim t^{-m}\|f\|_{p,m}\ee^{-(1-t^2)|z|^2/2}.
\]
For fixed $t<1$, the right-hand side is integrable after taking its $p$th power. If a bounded sequence tends to zero locally uniformly, so do its dilates; dominated convergence then gives convergence of the dilated sequence in $\F_m^p$. Compactness follows from Lemma~\ref{lem:compact}. To prove \eqref{eq:dilationlocal}, write
\begin{align*}
 ((\Id-R_t)f)^{(k)}(z)
 &=(1-t^k)f^{(k)}(z)
       +t^k\big(f^{(k)}(z)-f^{(k)}(tz)\big)\\
 &=(1-t^k)f^{(k)}(z)
       +t^k\int_t^1 zf^{(k+1)}(sz)\dd s.
\end{align*}
For $|z|\le R$ and $t\ge1/2$, both derivatives in the last line are bounded by \eqref{eq:point}; the integral has length $1-t$. Consequently the expression is at most $C_{R,k,p,m}(1-t)\|f\|_{p,m}$, uniformly on the unit ball of $\F_m^p$, proving \eqref{eq:dilationlocal}.
\end{proof}

\begin{lemma}\label{lem:compacttails}
Fix a quasi-normed space $X$ and an exponent $0<q<\infty$.
\begin{enumerate}
\item Let $K:X\to L^q(\C,d\tau)$ be compact, where $\tau$ is a positive measure. With $\chi_R=\mathbf1_{\{|z|>R\}}$, we have
\[
                    \|\chi_R K\|\longrightarrow0\quad(R\to\infty).
\]
\item If $K:X\to\F_\ell^q$ is compact, then
\begin{equation}\label{eq:compactpointtails}
 \sup_{\|f\|_X\le1}\sup_{|z|>R}
   |Kf(z)|(1+|z|)^\ell\ee^{-|z|^2/2}\longrightarrow0
   \quad(R\to\infty).
\end{equation}
\end{enumerate}
\end{lemma}
\begin{proof}
For fixed $h\in L^q(d\tau)$, dominated convergence shows $\|\chi_Rh\|_{L^q}\to0$. Compactness of $K$ provides a finite $\varepsilon$-net $h_1,\ldots,h_M$ for the image of the unit ball. Since multiplication by $\chi_R$ is contractive and the $\alpha$th power of the $L^q$ quasi-norm is subadditive for $\alpha=\min\{1,q\}$,
\[
 \|\chi_RK\|^\alpha
 \le\varepsilon^\alpha+
                  \max_{1\le i\le M}\|\chi_Rh_i\|_{L^q}^\alpha.
\]
Let $R\to\infty$ and then $\varepsilon\to0$ to obtain (a). For (b), apply \eqref{eq:local} with $k=0$ to $h\in\F_\ell^q$:
\[
 |h(z)|^q(1+|z|)^{\ell q}\ee^{-q|z|^2/2}
 \lesssim\int_{D(z,1)}|h(w)|^q(1+|w|)^{\ell q}
                         \ee^{-q|w|^2/2}\dd A(w).
\]
For $|z|>R$, the right-hand side is bounded by the weighted norm integral over $|w|>R-1$, which tends to zero as $R\to\infty$. Thus $(1+|z|)^\ell|h(z)|\ee^{-|z|^2/2}\to0$ as $|z|\to\infty$ for each fixed $h\in\F_\ell^q$. Choose a finite $\varepsilon$-net $h_1,\ldots,h_M$ for the relatively compact image of the unit ball under $K$ in $\F_\ell^q$. Uniform boundedness of point evaluations gives
\begin{align*}
 &\sup_{\|f\|_X\le1}\sup_{|z|>R}
                     |Kf(z)|(1+|z|)^\ell\ee^{-|z|^2/2}
 \le C\varepsilon+
       \max_{1\le i\le M}\sup_{|z|>R}
                            |h_i(z)|(1+|z|)^\ell\ee^{-|z|^2/2}.
\end{align*}
Each $h_i$ has the vanishing tail just established. Taking $R\to\infty$ and then $\varepsilon\to0$ completes (b).
\end{proof}

\section{Proofs of the Carleson Theorems}\label{sec:carleson}
The measure in \eqref{eq:jet} will be discretized on a lattice. We then use its disk masses to prove the boundedness and compactness criteria and the essential-norm estimates of Theorems~\ref{thm:carlesonup} and \ref{thm:carlesondown}.

\subsection{Disk Masses and Lattice Estimates}
Fix $\rho>0$. Let $\Lambda=\rho(\mathbb Z+i\mathbb Z)$ and let $Q_w$ be the half-open square of side length $\rho$ centered at $w\in\Lambda$. These squares partition $\C$ and satisfy $Q_w\subset D(w,\rho)$. The disks $D(w,3\rho)$ have bounded overlap. For a positive locally finite measure $\sigma$, write
\[
 b_w=\sigma(D(w,\rho)).
\]
Then, for $0<r<\infty$,
\begin{equation}\label{eq:lattice}
 \|B_\rho\sigma\|_{L^r}\asymp\|(b_w)\|_{\ell^r},
 \qquad
 \|B_\rho\sigma\|_\infty\asymp\sup_w b_w.
\end{equation}
The constants depend on $\rho$. The corresponding vanishing-at-infinity conditions are equivalent.

Here is a direct justification of the radius and lattice assertions. For fixed $R,\rho>0$, choose finitely many offsets $a_1,\ldots,a_M$ such that $D(0,R)\subset\bigcup_{j=1}^MD(a_j,\rho)$. Then
\[
 B_R\sigma(z)\le\sum_{j=1}^M B_\rho\sigma(z+a_j).
\]
Interchanging $R$ and $\rho$, and using translation invariance of area measure and the finite-sum inequality for $r$th powers, proves radius independence of the $L^r$, boundedness and vanishing conditions. This argument is valid also for $r<1$. If $z\in Q_w$, then $b_w\le B_{2\rho}\sigma(z)$. Conversely, $D(z,\rho)$ meets only a fixed number of squares, and its mass is at most the sum of the corresponding $b_v$. Integration on the disjoint squares proves both inequalities in \eqref{eq:lattice}, using radius independence for the first one. Only a fixed number of neighboring indices is involved, so the same argument proves the vanishing assertion.

\begin{lemma}\label{lem:upper}
For the measure $\nu$ in \eqref{eq:jet}, set $b_w=\nu(D(w,\rho))$. There exist nonnegative quantities $A_w(f)$ satisfying
\begin{equation}\label{eq:discreteupper}
 \begin{aligned}
 \sum_w A_w(f)^p
   &\lesssim\|f\|_{p,m}^p,
 \qquad \|J_{\mathbf u,\mu}f\|_{L^q}^q
   &\lesssim\sum_w A_w(f)^q b_w.
 \end{aligned}
\end{equation}
For instance, take
\[
 A_w(f)^p=\int_{D(w,3\rho)}|f(z)|^p(1+|z|)^{mp}
                              \ee^{-p|z|^2/2}\dd A(z).
\]
\end{lemma}
\begin{proof}
Put $F(z)=|f(z)|^p(1+|z|)^{mp}\ee^{-p|z|^2/2}$. By Tonelli's theorem and the bounded overlap of the disks $D(w,3\rho)$,
\[
 \sum_w A_w(f)^p
 =\int_{\C}F(z)\sum_w\mathbf1_{D(w,3\rho)}(z)\dd A(z)
 \lesssim\|f\|_{p,m}^p.
\]
If $z\in Q_w$, then $|z-w|<\rho$ and $D(z,\rho)\subset D(w,3\rho)$. Apply Lemma~\ref{lem:local} on the smaller disk to obtain
\[
           (1+|z|)^m|Q_kf(z)|\lesssim A_w(f),\qquad 0\le k\le N.
\]
Substitute these bounds in \eqref{eq:algebra} and use the finite-sum inequality for the $q$th power. For $z\in Q_w$ this gives
\[
 |\du f(z)|^q\ee^{-q|z|^2/2}
 \lesssim A_w(f)^q
 \frac{\sum_{k=0}^N|P_{\mathbf u(z)}^{(k)}(\bar z)|^q}
      {(1+|z|)^{mq}},\qquad z\in Q_w.
\]
Integration over $Q_w$ with respect to $\mu$ gives
\begin{align*}
 \int_{Q_w}|\du f(z)|^q\ee^{-q|z|^2/2}\dd\mu(z)
 &\lesssim A_w(f)^q
   \int_{Q_w}\frac{\sum_{k=0}^N|P_{\mathbf u(z)}^{(k)}(\bar z)|^q}
                       {(1+|z|)^{mq}}\dd\mu(z)\\
 &\le A_w(f)^q b_w.
\end{align*}
The squares partition $\C$, so summing over $w$ gives the second inequality in \eqref{eq:discreteupper}.
\end{proof}

\subsection{\texorpdfstring{Proof of Theorem~\ref{thm:carlesonup}}{The Carleson theorem for p <= q}}
\begin{proof}
Suppose first that (\ref{cup:mass}) holds. Set $J=J_{\mathbf u,\mu}$ and $b_w=\nu(D(w,\rho))$ for $w\in\Lambda$. The discrete estimate \eqref{eq:discreteupper}, followed by $\ell^p\subset\ell^q$ for $p\le q$, gives
\begin{align*}
 \|Jf\|_{L^q}^q
 &\lesssim\sum_{w\in\Lambda}A_w(f)^q b_w\\
 &\le\left(\sup_{w\in\Lambda}b_w\right)
                       \sum_{w\in\Lambda}A_w(f)^q\\
 &\le\left(\sup_{w\in\Lambda}b_w\right)
                  \left(\sum_{w\in\Lambda}A_w(f)^p\right)^{q/p}\\
 &\lesssim\|B_\rho\nu\|_\infty\|f\|_{p,m}^q.
\end{align*}
Thus (\ref{cup:mu}) holds with the stated upper bound for the norm. The sequence inequality uses only $q/p\ge1$, so the argument applies without assuming either $p\ge1$ or $q\ge1$.

For the converse, assume (\ref{cup:mu}) and fix $w\in\C$. Apply the boundedness estimate for $J_{\mathbf u,\mu}$ to each $E_{w,k}$, $0\le k\le N$. Lemma~\ref{lem:atoms} bounds $\|E_{w,k}\|_{p,m}$ uniformly, while the local test comparison in Lemma~\ref{lem:algebra} yields
\begin{align}
 B_\rho\nu(w)
 &=\int_{D(w,\rho)}
     \frac{\sum_{k=0}^N|P_{\mathbf u(z)}^{(k)}(\bar z)|^q}
                         {(1+|z|)^{mq}}\dd\mu(z)\notag\\
 &\lesssim\sum_{k=0}^N\int_{D(w,\rho)}
         |\du E_{w,k}(z)|^q\ee^{-q|z|^2/2}\dd\mu(z)\notag\\
 &\le\sum_{k=0}^N\|JE_{w,k}\|_{L^q}^q
 \le\|J\|^q\sum_{k=0}^N\|E_{w,k}\|_{p,m}^q
 \lesssim\|J\|^q.
 \label{eq:necessityup}
\end{align}
Taking the supremum over $w$ proves (\ref{cup:mass}) and the reverse norm estimate. Hence (\ref{cup:mu}) and (\ref{cup:mass}) are equivalent. With $m=N=0$, $u_0=1$, and $\mu$ replaced by $\nu$, the same proof gives the $(p,q)$-Fock Carleson criterion, namely (\ref{cup:nu})$\Leftrightarrow$(\ref{cup:mass}). No separate scalar result is required.

Turn now to compactness. Assume (\ref{cup:vanmass}). Since $\nu$ is locally finite, for every fixed $R$ we have
\[
 \sup_{|w|\le R}B_\rho\nu(w)
       \le\nu(D(0,R+\rho))<\infty.
\]
The assumed vanishing at infinity therefore makes $B_\rho\nu$ bounded on all of $\C$, and the boundedness criterion applies to $J$. Let $(f_j)$ be bounded in $\F_m^p$ and converge to zero uniformly on compact subsets of $\C$, and put $M=\sup_j\|f_j\|_{p,m}$. Split the sum in Lemma~\ref{lem:upper} at radius $R$:
\begin{align*}
 \|Jf_j\|_{L^q}^q
 &\lesssim\sum_{|w|\le R}A_w(f_j)^q b_w
          +\sum_{|w|>R}A_w(f_j)^q b_w\\
 &\lesssim\sum_{|w|\le R}A_w(f_j)^q b_w
          +M^q\sup_{|w|>R}b_w.
\end{align*}
For fixed $R$, only finitely many lattice points appear in the first sum, and local uniform convergence implies $A_w(f_j)\to0$ for each of them. Given $\varepsilon>0$, choose $R$ with $\sup_{|w|>R}b_w<\varepsilon$ and then let $j\to\infty$. The displayed bound gives $\limsup_j\|Jf_j\|_{L^q}^q\lesssim M^q\varepsilon$. Letting $\varepsilon\to0$ proves (\ref{cup:vanmu}).

Conversely, assume (\ref{cup:vanmu}) and choose any sequence $|w_j|\to\infty$. For each $k$, the functions $E_{w_j,k}$ are uniformly bounded in $\F_m^p$ and converge to zero uniformly on compact subsets of $\C$ by Lemma~\ref{lem:atoms}. Thus
\[
 \nu(D(w_j,\rho))\lesssim
               \sum_{k=0}^N\|JE_{w_j,k}\|_{L^q}^q\longrightarrow0,
\]
where we used the first two inequalities of \eqref{eq:necessityup}. Since the sequence of centers was arbitrary, (\ref{cup:vanmass}) follows. Repeating both implications with $m=N=0$, $u_0=1$, and measure $\nu$ proves (\ref{cup:vannu})$\Leftrightarrow$(\ref{cup:vanmass}). Lemma~\ref{lem:compact} then identifies the sequential vanishing condition with compactness of $J$.

Assume for the essential-norm estimate that $J$ is bounded, and write
\[
 Y=L^q(\C,\ee^{-q|z|^2/2}d\mu),\qquad
 L_\rho=\limsup_{|w|\to\infty}B_\rho\nu(w),\qquad
 \chi_R=\mathbf1_{\{|z|>R\}}.
\]
If $Q_w$ meets $\{|z|>R\}$, its center satisfies $|w|>R-\rho$. Restrict the proof of Lemma~\ref{lem:upper} to these portions of the squares to obtain
\begin{align}
 \|\chi_RJf\|_Y^q
 &\lesssim\sum_{\substack{w\in\Lambda\\|w|>R-\rho}}
                                      A_w(f)^q b_w\notag\\
 &\lesssim\left(\sup_{\substack{w\in\Lambda\\|w|>R-\rho}}b_w\right)
                                      \|f\|_{p,m}^q.
 \label{eq:embeddingtail}
\end{align}
For each $R$, the truncated operator $J_R=(1-\chi_R)J$ is compact. Indeed, every bounded sequence in $\F_m^p$ has, by Montel's theorem, a subsequence converging uniformly on compact subsets of $\C$. Its derivatives converge uniformly on $\overline{D(0,R)}$; the coefficients $u_j$ are bounded there and $\mu(\overline{D(0,R)})<\infty$. The corresponding images under $J_R$ therefore converge in $Y$. Together with \eqref{eq:embeddingtail}, this gives
\[
 \|J\|_{\mathrm e}^q
 \le\liminf_{R\to\infty}\|J-J_R\|^q\lesssim L_\rho.
\]
The lattice supremum is no larger than the supremum over all centers, so the quantity on the right is the $L_\rho$ appearing in the statement.

For the lower estimate, fix any compact $K:\F_m^p\to Y$ and take $|w|>R+\rho$. The test disk $D(w,\rho)$ lies entirely in $\{|z|>R\}$, so
\begin{align*}
 B_\rho\nu(w)
 &\lesssim\sum_{k=0}^N\|\chi_RJE_{w,k}\|_Y^q\\
 &\lesssim\sum_{k=0}^N\left(
       \|(J-K)E_{w,k}\|_Y^q+\|\chi_RKE_{w,k}\|_Y^q\right)\\
 &\lesssim\|J-K\|^q+\|\chi_RK\|^q.
\end{align*}
The constants do not depend on $R$, $w$, or $K$. Take the upper limit over $|w|\to\infty$ before sending $R\to\infty$. Lemma~\ref{lem:compacttails}(a) then removes the term $\|\chi_RK\|^q$, leaving $L_\rho\lesssim\|J-K\|^q$. The infimum over compact $K$ gives the lower comparison in \eqref{eq:essentialmeasure}. The uniform tail estimate is available even for $p<1$.

For the essential-norm estimates of integral operators, we also need compact approximants whose defining expressions preserve entire functions. The operators $JR_t$ have this property and are compact by Lemma~\ref{lem:dilations}. For $h=(\Id-R_t)f$, $\|f\|_{p,m}\le1$, and $1/2\le t<1$, uniform boundedness of $R_t$ gives $\|h\|_{p,m}\le C$. On the fixed disk $|z|\le R$, boundedness of the coefficients and local finiteness of $\mu$ yield
\[
 \|(1-\chi_R)Jh\|_Y^q
 \le C_R\sum_{k=0}^N\sup_{|z|\le R}|h^{(k)}(z)|^q.
\]
The supremum of this local bound over $\|f\|_{p,m}\le1$ tends to zero as $t\to1$ by \eqref{eq:dilationlocal}. On the complement of the disk, apply \eqref{eq:embeddingtail}. The two regions partition the integral defining $\|Jh\|_Y^q$, and hence
\[
 \limsup_{t\to1}\|J(\Id-R_t)\|^q
       \lesssim\sup_{\substack{w\in\Lambda\\|w|>R-\rho}}b_w.
\]
Now send $R\to\infty$ to obtain
\begin{equation}\label{eq:holomorphicapprox}
 \limsup_{t\to1}\|J(\Id-R_t)\|^q\lesssim L_\rho.
\end{equation}
In this argument $t\to1$ is taken for each fixed $R$; only then is $R$ allowed to tend to infinity.
\end{proof}

\subsection{\texorpdfstring{Proof of Theorem~\ref{thm:carlesondown}}{The Carleson theorem for q < p}}
\begin{proof}
Suppose (\ref{cdown:mass}) holds and put $r=\frac{p}{p-q}>1$. Apply the lattice comparison and H\"older's inequality with conjugate exponents $p/q$ and $r$ to obtain
\begin{align}
 \|J_{\mathbf u,\mu}f\|_{L^q}^q
 &\lesssim\sum_wA_w(f)^q b_w\notag\\
 &\le\left(\sum_wA_w(f)^p\right)^{q/p}
             \left(\sum_wb_w^r\right)^{1/r}\notag\\
 &\lesssim\|f\|_{p,m}^q\|(b_w)\|_{\ell^r}.
 \label{eq:downupper}
\end{align}
For necessity of (\ref{cdown:mass}), testing one atom at a time as in \eqref{eq:necessityup} controls only $\sup_w b_w$. To obtain $\ell^r$ summability, we test randomized finite linear combinations of the atoms.

Take a finitely supported sequence $c=(c_w)$ and independent Rademacher functions $\varepsilon_w(t)$ on $[0,1]$. For each $k\in\{0,\ldots,N\}$, define
\[
                F_{t,k}=\sum_{w\in\Lambda}
                           c_w\varepsilon_w(t)E_{w,k}.
\]
The synthesis estimate in Lemma~\ref{lem:atoms} is uniform in the signs and therefore gives
\[
                  \|F_{t,k}\|_{p,m}\lesssim\|c\|_{\ell^p}.
\]
Now assume $J_{\mathbf u,\mu}$ is bounded. We use the Khintchine inequality, valid for every $q>0$ and every finite complex scalar family:
\[
 \int_0^1\left|\sum_w\varepsilon_w(t)\alpha_w\right|^q\dd t
 \asymp_q\left(\sum_w|\alpha_w|^2\right)^{q/2}.
\]
Apply this inequality at each $z$ with $\alpha_w=c_w\du E_{w,k}(z)$, and integrate the operator bound for $F_{t,k}$ over $t$. Tonelli's theorem then gives
\begin{align}
 &\int_{\C}\left(\sum_w|c_w\du E_{w,k}(z)|^2\right)^{q/2}
                              \ee^{-q|z|^2/2}\dd\mu(z)\notag\\
 &\qquad\asymp\int_0^1\int_{\C}
          |\du F_{t,k}(z)|^q\ee^{-q|z|^2/2}\dd\mu(z)\dd t\notag\\
 &\qquad\lesssim\|J_{\mathbf u,\mu}\|^q\|c\|_{\ell^p}^q.
 \label{eq:random}
\end{align}
All sums in this argument are finite, and the integrands are nonnegative, so Tonelli's theorem requires no preliminary integrability assumption.

The disks $D(w,\rho)$ have a fixed overlap bound $M_\rho$. At a given point $z$, at most $M_\rho$ of their indicators are nonzero. For $0<q<2$, H\"older's inequality on this finite set gives the following estimate; for $q\ge2$, it is the inclusion $\ell^2\subset\ell^q$: for arbitrary nonnegative numbers $a_w$,
\[
 \sum_w a_w^q\mathbf1_{D(w,\rho)}(z)
 \le M_\rho^{\max\{1-q/2,0\}}
       \left(\sum_w a_w^2\mathbf1_{D(w,\rho)}(z)\right)^{q/2}.
\]
Choose $a_w=|c_w\du E_{w,k}(z)|$, multiply by $\ee^{-q|z|^2/2}$, and integrate against $\mu$. Removing the indicators from the square sum can only enlarge it, so \eqref{eq:random} implies
\begin{align*}
 &\sum_w|c_w|^q\int_{D(w,\rho)}
           |\du E_{w,k}(z)|^q\ee^{-q|z|^2/2}\dd\mu(z)
\lesssim
              \|J_{\mathbf u,\mu}\|^q\|c\|_{\ell^p}^q.
\end{align*}
Summing over $k$ and applying the local equivalence \eqref{eq:localtests} to each disk gives
\begin{equation}\label{eq:sequencefunctional}
        \sum_w|c_w|^q b_w
           \lesssim\|J_{\mathbf u,\mu}\|^q
                        \left(\sum_w|c_w|^p\right)^{q/p}.
\end{equation}
To identify the sequence norm, fix a finite $E\subset\Lambda$ and choose
\[
         |c_w|^q=b_w^{r-1}\mathbf1_E(w).
\]
The identity $(r-1)p/q=r$ turns \eqref{eq:sequencefunctional} into
\[
 \sum_{w\in E}b_w^r
 \lesssim\|J_{\mathbf u,\mu}\|^q
                \left(\sum_{w\in E}b_w^r\right)^{q/p}.
\]
For a nonzero sum, divide by its $q/p$ power and use $1-q/p=1/r$ to obtain
\[
 \left(\sum_{w\in E}b_w^r\right)^{1/r}
 \lesssim\|J_{\mathbf u,\mu}\|^q.
\]
The same inequality is immediate when the sum vanishes. Letting finite sets $E$ increase to $\Lambda$ and using monotone convergence yields
\begin{equation}\label{eq:downlower}
          \|(b_w)\|_{\ell^r}\lesssim\|J_{\mathbf u,\mu}\|^q.
\end{equation}
Combining \eqref{eq:lattice}, \eqref{eq:downupper}, and \eqref{eq:downlower} proves (\ref{cdown:mu})$\Leftrightarrow$(\ref{cdown:mass}) and the claimed norm equivalence.

Next let $(f_j)$ be bounded in $\F_m^p$ and converge to zero uniformly on compact subsets of $\C$. Under (\ref{cdown:mu}), the disk masses lie in $\ell^r$ by the criterion just proved. Splitting the lattice sum into a finite part and a tail gives
\begin{align*}
 \|J_{\mathbf u,\mu}f_j\|_{L^q}^q
 &\lesssim\sum_{|w|\le R}A_w(f_j)^q b_w\\
 &\quad+\left(\sum_{|w|>R}A_w(f_j)^p\right)^{q/p}
                \left(\sum_{|w|>R}b_w^r\right)^{1/r}\\
 &\lesssim\sum_{|w|\le R}A_w(f_j)^q b_w
       +\sup_j\|f_j\|_{p,m}^q
                \left(\sum_{|w|>R}b_w^r\right)^{1/r}.
\end{align*}
Given $\varepsilon>0$, choose $R$ so that $(\sum_{|w|>R}b_w^r)^{1/r}<\varepsilon$. For this fixed $R$, the first sum has finitely many terms, each tending to zero because $f_j\to0$ locally uniformly. Taking the upper limit as $j\to\infty$ yields
\[
 \limsup_{j\to\infty}\|J_{\mathbf u,\mu}f_j\|_{L^q}^q
       \lesssim\varepsilon\sup_j\|f_j\|_{p,m}^q.
\]
Let $\varepsilon\to0$. This proves (\ref{cdown:vanmu}), which implies compactness by Lemma~\ref{lem:compact}. Conversely, (\ref{cdown:vanmu}) includes boundedness as part of the definition of a vanishing Carleson measure.

Finally, set $m=N=0$ and $u_0=1$, replacing $\mu$ by $\nu$ in the criterion above. This gives (\ref{cdown:nu})$\Leftrightarrow$(\ref{cdown:mass}). The scalar compactness argument also identifies (\ref{cdown:nu}) with (\ref{cdown:vannu}), completing the equivalence of the five conditions.
\end{proof}

\section{Proofs of the Operator Theorems}\label{sec:operators}
We apply the Carleson criteria to the operators introduced above. We begin with generalized integral operators and then identify the coefficients that determine their essential norms and Fredholm indices. The remaining subsections address differential equations and composition--differentiation operators.

\subsection{Generalized Integral Operators}\label{subsec:integralproof}
For a generalized integral operator, write
\begin{equation}\label{eq:Ag}
 A_{\mathbf g}(z)=\frac{S_{\mathbf g}(z)}{(1+|z|)^d},
 \qquad d=n+m-\ell.
\end{equation}
To prove Theorem~\ref{thm:integral}, we first establish degree criteria for
the coefficients in terms of $A_{\mathbf g}$. The following lemma
extends the polynomial argument of
\cite[Section~5]{AHTYY} to the present weights. Retaining every derivative
of the coefficient polynomial also covers negative weight exponents and
the strict integrability threshold.

\begin{lemma}\label{lem:rigidity}
Let $d\in\mathbb R$, and define $A_{\mathbf g}$ by \eqref{eq:Ag}. The following assertions hold:
\begin{enumerate}
\item $A_{\mathbf g}$ is bounded if and only if each nonzero $g_j$ is a polynomial with $\deg g_j+j\le d$;
\item $A_{\mathbf g}(z)\to0$ as $|z|\to\infty$ if and only if each nonzero $g_j$ is a polynomial with $\deg g_j+j<d$;
\item For $0<s<\infty$, $A_{\mathbf g}\in L^s(dA)$ if and only if each nonzero $g_j$ is a polynomial with $\deg g_j+j<d-\frac2s$.
\end{enumerate}
\end{lemma}
\begin{proof}
We begin with (a). Assume that $M=\|A_{\mathbf g}\|_\infty<\infty$.
For $0\le j\le N$, Taylor's formula in the polynomial variable gives
\begin{align*}
 g_j(z)
 &=\frac{1}{j!}P_{\mathbf g(z)}^{(j)}(0)
 =\frac1{j!}\sum_{k=j}^N
       \frac{P_{\mathbf g(z)}^{(k)}(\bar z)}{(k-j)!}(-\bar z)^{k-j}.
\end{align*}
Consequently,
\[
 |g_j(z)|\lesssim (1+|z|)^{N-j}S_{\mathbf g}(z)
                 \le M(1+|z|)^{d+N-j}.
\]
Write $g_j(z)=\sum_{h=0}^\infty b_{j,h}z^h$. Cauchy's coefficient estimate on $|z|=R$, with $R\ge1$, yields
\[
 |b_{j,h}|\le R^{-h}\max_{|z|=R}|g_j(z)|
            \lesssim M R^{d+N-j-h}.
\]
If $h>d+N-j$, letting $R\to\infty$ gives $b_{j,h}=0$. Thus every $g_j$ is a polynomial; when $d+N-j<0$, all of its coefficients vanish. This preliminary estimate proves polynomiality, but we still have to obtain the sharper degree bound $\deg g_j+j\le d$.

If all $g_j$ vanish, there is nothing more to prove. Otherwise, put
\[
 D=\max_{g_j\ne0}(\deg g_j+j),\qquad
 F(z)=P_{\mathbf g(z)}(\bar z)=\sum_{j=0}^Ng_j(z)\bar z^j.
\]
Let $c_j$ be the leading coefficient of $g_j$ when $\deg g_j+j=D$. The terms of total degree $D$ give
\begin{equation}\label{eq:leading}
 \begin{aligned}
 R^{-D}F(R\ee^{i\vartheta})
   &=H_D(\vartheta)+O(R^{-1}),
 \quad H_D(\vartheta)
   &=\sum_{\deg g_j+j=D}c_j\ee^{i(D-2j)\vartheta},
 \end{aligned}
\end{equation}
uniformly in $\vartheta\in[0,2\pi]$; when $D=0$, there is no error term.
Since the frequencies $D-2j$ are distinct, orthogonality yields
\[
 \frac1{2\pi}\int_0^{2\pi}|H_D(\vartheta)|^2\dd\vartheta
          =\sum_{\deg g_j+j=D}|c_j|^2>0.
\]
Choose $\vartheta_0$ with $H_D(\vartheta_0)\ne0$. By \eqref{eq:leading}, for all sufficiently large $R$,
\[
 \tfrac12|H_D(\vartheta_0)|R^D
       \le|F(R\ee^{i\vartheta_0})|
       \le S_{\mathbf g}(R\ee^{i\vartheta_0})
       \le M(1+R)^d.
\]
This is impossible if $D>d$. Therefore $D\le d$, as required. In particular, boundedness forces $\mathbf g=0$ if $d<0$.

Conversely, suppose the nonzero $g_j$ are polynomials with
$\deg g_j+j\le d$. Retain the notation $D$ from above. Then
\[
 P_{\mathbf g(z)}^{(k)}(\bar z)
   =\sum_{j=k}^N\sum_{h=0}^{\deg g_j}
                  \frac{j!}{(j-k)!}b_{j,h}z^h\bar z^{j-k}.
\]
Each nonzero monomial has total degree $h+j-k\le D-k$. It follows that
\[
 S_{\mathbf g}(z)\le C_{\mathbf g}(1+|z|)^D,
 \qquad A_{\mathbf g}(z)\le C_{\mathbf g}(1+|z|)^{D-d}.
\]
The latter expression is bounded because $D\le d$. This proves (a). The constant $C_{\mathbf g}$ here may depend on the fixed polynomial coefficients; it is used only to establish finiteness and decay.

For (b), assume that $A_{\mathbf g}(z)\to0$. This continuous function is then bounded on $\C$, so (a) gives polynomial coefficients and $D\le d$. If $D=d$, the preceding lower estimate along $R\ee^{i\vartheta_0}$ gives
\[
 A_{\mathbf g}(R\ee^{i\vartheta_0})
 \ge\frac{|F(R\ee^{i\vartheta_0})|}{(1+R)^d}
 \ge\tfrac12|H_D(\vartheta_0)|\frac{R^d}{(1+R)^d}.
\]
The right-hand side has a positive limit, contrary to the assumed decay. Thus $D<d$. If instead $D<d$, the polynomial upper bound immediately gives $A_{\mathbf g}(z)\to0$.

For (c), let $A_{\mathbf g}\in L^s(dA)$. In Lemma~\ref{lem:submeanjet}, choose $W(z)=(1+|z|)^{-d}$,
$\eta(z)=z$, $\beta=0$, $t=s$, and radius one. We obtain
\[
 A_{\mathbf g}(w)^s
       \lesssim\int_{D(w,1)}A_{\mathbf g}(z)^s\dd A(z)
       \le\|A_{\mathbf g}\|_{L^s}^s.
\]
Thus $A_{\mathbf g}$ is bounded, and (a) shows that the coefficients
are polynomials. Suppose $\mathbf g\ne0$. By continuity, the nonzero
value $H_D(\vartheta_0)$ persists on an interval $J$ of positive length:
$|H_D|\ge2c$ there for some $c>0$. The uniform expansion
\eqref{eq:leading} then provides $R_0>1$ such that
\[
 |F(R\ee^{i\vartheta})|\ge cR^D
       \quad(R\ge R_0,\ \vartheta\in J).
\]
Polar integration now gives
\begin{align}
 \int_{\C}A_{\mathbf g}(z)^s\dd A(z)
 &\ge\int_J\int_{R_0}^\infty
       \frac{|F(R\ee^{i\vartheta})|^s}{(1+R)^{ds}}R\dd R\dd\vartheta
       \notag\\
 &\ge c'\int_{R_0}^\infty R^{s(D-d)+1}\dd R.
 \label{eq:radialthreshold}
\end{align}
Finiteness of the last integral requires
$s(D-d)+1<-1$, or $D<d-\frac2s$; equality produces the divergent
$R^{-1}$ integral. Conversely, under this strict inequality, the
polynomial upper bound yields
\[
 \int_{\C}A_{\mathbf g}(z)^s\dd A(z)
 \le C_{\mathbf g}^s\,2\pi\int_0^\infty
                      (1+R)^{s(D-d)}R\dd R<\infty.
\]
This completes (c). The zero coefficient vector satisfies each assertion directly.
\end{proof}

\begin{proof}[Proof of Theorem \ref{thm:integral}]
Apply \eqref{eq:primitive} to
$H=\sum_{j=0}^Ng_jf^{(j)}$. This gives
\begin{equation}\label{eq:integralnorm}
 \|\Tg f\|_{q,\ell}^q\asymp
 \int_{\C}\left|\sum_{j=0}^Ng_j(z)f^{(j)}(z)\right|^q
          (1+|z|)^{(\ell-n)q}\ee^{-q|z|^2/2}\dd A(z).
\end{equation}
Set $d\mu(z)=(1+|z|)^{(\ell-n)q}dA(z)$ and
$J=J_{\mathbf g,\mu}$. By \eqref{eq:integralnorm}, the boundedness of
$T=\Tg$ is equivalent to the $(\mathbf g,p,q;m)$-Carleson property of
$\mu$; moreover, $\|T\|\asymp\|J\|$ whenever these conditions hold.
The associated measure $\nu$ in \eqref{eq:jet} is then
\begin{align*}
 d\nu(z)
 &=\frac{\sum_{k=0}^N|P_{\mathbf g(z)}^{(k)}(\bar z)|^q}
                         {(1+|z|)^{mq}}
                            (1+|z|)^{(\ell-n)q}\dd A(z)\\
 &=\frac{\sum_{k=0}^N|P_{\mathbf g(z)}^{(k)}(\bar z)|^q}
                         {(1+|z|)^{dq}}\dd A(z)
 \asymp A_{\mathbf g}(z)^q\dd A(z).
\end{align*}
Here $d=n+m-\ell$, so the two Sobolev weights combine into the
exponent $-dq$.

Suppose first that $p\le q$. Theorem~\ref{thm:carlesonup} reduces
boundedness to
\[
 \sup_{w\in\C}\int_{D(w,1)}A_{\mathbf g}(z)^q\dd A(z)<\infty.
\]
To compare the local disk integrals with $A_{\mathbf g}$ itself, fix $w$ and
$0\le k\le N$. The map
$z\mapsto P_{\mathbf g(z)}^{(k)}(\bar w)$ is entire. Applying the
submean inequality to this map and then expanding in the polynomial
variable at $\bar z$, we obtain
\begin{align*}
 |P_{\mathbf g(w)}^{(k)}(\bar w)|^q
 &\lesssim\int_{D(w,1)}
                   |P_{\mathbf g(z)}^{(k)}(\bar w)|^q\dd A(z)\\
 &=\int_{D(w,1)}
   \left|\sum_{v=k}^N\frac{P_{\mathbf g(z)}^{(v)}(\bar z)}{(v-k)!}
                                 (\bar w-\bar z)^{v-k}\right|^q\dd A(z)\\
 &\lesssim\int_{D(w,1)}
                \sum_{v=k}^N|P_{\mathbf g(z)}^{(v)}(\bar z)|^q\dd A(z).
\end{align*}
On $D(w,1)$ we have $1+|z|\asymp1+|w|$. Divide the preceding
estimates by $(1+|w|)^{dq}$ and sum over $k$ to get
\begin{equation}\label{eq:Agdisk}
 A_{\mathbf g}(w)^q
 \lesssim\int_{D(w,1)}A_{\mathbf g}(z)^q\dd A(z)
 \le\pi\sup_{z\in D(w,1)}A_{\mathbf g}(z)^q.
\end{equation}
Alternatively, the first inequality follows from
Lemma~\ref{lem:submeanjet} with $\eta(z)=z$ and $\beta=0$.
Equation~\eqref{eq:Agdisk} therefore gives
\[
 \|B_1\nu\|_\infty\asymp\|A_{\mathbf g}\|_\infty^q.
\]
Thus $T$ is bounded exactly when $A_{\mathbf g}$ is bounded.
Lemma~\ref{lem:rigidity}(a) translates this into
$\deg g_j+j\le d$ for every nonzero $g_j$, and the same estimates give
$\|T\|\asymp\|A_{\mathbf g}\|_\infty$.

For compactness, note that the expressions defining $T$ and $J$
preserve local uniform convergence of entire functions. Indeed,
derivatives converge locally uniformly, and the integration paths
defining $T$ remain in a compact set when their endpoints do. Hence, for any
bounded sequence $(f_j)$ tending to zero locally uniformly,
\eqref{eq:integralnorm} implies
\[
 \|Tf_j\|_{q,\ell}\to0
            \quad\Longleftrightarrow\quad\|Jf_j\|_{L^q}\to0.
\]
Lemma~\ref{lem:compact} therefore makes compactness of $T$ equivalent
to compactness of $J$. Theorem~\ref{thm:carlesonup} characterizes the
latter by $B_1\nu(w)\to0$. By \eqref{eq:Agdisk}, this is equivalent to
$A_{\mathbf g}(w)\to0$: in the reverse direction, the disk $D(w,1)$
escapes every fixed compact set as $|w|\to\infty$. The strict degree
conditions in part (a) now follow from Lemma~\ref{lem:rigidity}(b).

Now turn to $q<p$. Set $r=\frac{p}{p-q}>1$ and $s=qr=\frac{pq}{p-q}$. By Theorem~\ref{thm:carlesondown}, boundedness and compactness of $J$ are equivalent, and they hold precisely when
\[
 G_1(w):=\int_{D(w,1)}A_{\mathbf g}(z)^q\dd A(z)
                                      \quad\text{belongs to }L^r.
\]
The same is true of $T$ by \eqref{eq:integralnorm} and the compactness argument just given. The lower estimate in \eqref{eq:Agdisk} implies
\[
 \int_{\C}A_{\mathbf g}(w)^{qr}\dd A(w)
                    \lesssim\int_{\C}G_1(w)^r\dd A(w).
\]
For the reverse estimate, H\"older's inequality on the disk and Tonelli's theorem yield
\begin{align*}
 \int_{\C}G_1(w)^r\dd A(w)
 &\le\pi^{r-1}\int_{\C}\int_{D(w,1)}
                       A_{\mathbf g}(z)^{qr}\dd A(z)\dd A(w)\\
 &=\pi^{r-1}\int_{\C}A_{\mathbf g}(z)^{qr}
                         \left(\int_{D(z,1)}\dd A(w)\right)\dd A(z)\\
 &=\pi^r\int_{\C}A_{\mathbf g}(z)^{qr}\dd A(z).
\end{align*}
Consequently,
\begin{equation}\label{eq:averageLr}
 \|G_1\|_{L^r}\asymp\|A_{\mathbf g}\|_{L^s}^{q}.
\end{equation}
By Lemma~\ref{lem:rigidity}(c), this is finite if and only if
\[
 \deg g_j+j<d-\frac2s
            =d-2\left(\frac1q-\frac1p\right),\qquad 0\le j\le N.
\]
This gives part (b), with a strict inequality at the threshold. The norm
statement follows from Theorem~\ref{thm:carlesondown} and
\eqref{eq:averageLr}:
\[
 \|T\|^q\asymp\|J\|^q\asymp\|B_1\nu\|_{L^r}
                       \asymp\|A_{\mathbf g}\|_{L^s}^q.
\]
Combining the two exponent ranges, we obtain
\[
 \|\Tg\|\asymp
 \begin{cases}
 \|A_{\mathbf g}\|_\infty,&p\le q,\\
 \|A_{\mathbf g}\|_{L^s},&q<p,
 \end{cases}
\]
where the implied constants do not depend on $\mathbf g$.
\end{proof}

\begin{proof}[Proof of Corollary \ref{cor:chalmoukis}]
Set $g_j=a_jg^{(n-j)}$. Suppose first that $T_{g,\mathbf a}=0$. For every $t\in\C$, the function $f_t(z)=\ee^{tz}$ belongs to the source space, because its linear exponential growth is dominated by the Gaussian in \eqref{eq:norm}. Differentiating $T_{g,\mathbf a}f_t=0$ exactly $n$ times gives
\[
 0=\sum_{j=0}^{n-1}g_j(z)t^j\ee^{tz},
 \qquad\text{hence}\qquad\sum_{j=0}^{n-1}g_j(z)t^j=0.
\]
For each fixed $z$, the polynomial in $t$ vanishes identically;
hence $g_j(z)=0$ for all $j$. In particular,
$a_{j_*}\ne0$ and $n-j_*=h$ imply $g^{(h)}=0$, or $\deg g<h$.
Conversely, if $\deg g<h$, then $n-j\ge h$ for every index $j$ with $a_j\ne0$. All corresponding derivatives of $g$ vanish, so
$T_{g,\mathbf a}=0$.

Suppose next that $T_{g,\mathbf a}$ is bounded and nonzero. By
Theorem~\ref{thm:integral}, $g^{(h)}$ is a polynomial, and therefore so is
$g$. Set $D=\deg g\ge h$. For each index $j$ with $a_j\ne0$ and $h_j=n-j\ge h$,
\[
 \deg(a_jg^{(h_j)})+j=D+n-2h_j\le D+n-2h,
\]
provided that the derivative is nonzero. Equality holds for $j=j_*$, so the three degree conditions of Theorem~\ref{thm:integral}
reduce to the asserted criteria.
\end{proof}

For the usual Volterra operator $T_gf=I(g'f)$, a nonconstant polynomial symbol of degree $D$ gives the thresholds $D\le2+m-\ell$ for boundedness when $p\le q$, and $D<2+m-\ell-2\big(\frac1q-\frac1p\big)$ when $q<p$. Compactness in the first range requires strict inequality. The companion $J_gf=I(gf')$ has only constant bounded symbols when $m=\ell$ and $p\le q$. A positive Sobolev gap also permits compact companions; for example, $J_1:\F_1^p\to\F_0^p$ is compact.

The restriction on $\mathbf a$ in Corollary~\ref{cor:chalmoukis} cannot be omitted. Take $n=2$, $\mathbf a=(1,0)$ and $g(z)=z^4$. Then
\[
              T_{g,\mathbf a}f=I^2(12z^2f)
\]
is bounded from $\F_m^p$ to $\F_m^q$ when $p\le q$, although $g$ is not quadratic. Here $h=2$, and the correct bound is $\deg g\le4$.

\subsection{Essential Norms and Fredholm Properties}\label{subsec:essentialproof}
We prove Theorems~\ref{thm:essentialintegral} and
\ref{thm:fredholmsymbol} by comparing the essential norms of the differential operator and its $n$-fold integral. A calculation in normalized monomial bases then gives the Hilbert-space formula and the Fredholm criterion.

\begin{proof}[Proof of Theorem \ref{thm:essentialintegral}: essential-norm comparison]
Write $T=\Tg$, $J=J_{\mathbf g,\mu}$, and $d\mu=(1+|z|)^{(\ell-n)q}dA$, as in \eqref{eq:integralnorm}. Then $J=D^nT$, with the derivative regarded as a map into $Y=L^q(\C,\ee^{-q|z|^2/2}d\mu)$. This derivative is bounded on the whole target $\F_\ell^q$. Indeed, define
\[
 \Pi_{n-1}h=\sum_{j=0}^{n-1}\frac{h^{(j)}(0)}{j!}z^j.
\]
The point-evaluation estimates show that $\Pi_{n-1}$ is bounded, while
$I^nh^{(n)}=h-\Pi_{n-1}h$. Lemma~\ref{lem:LP} therefore yields
\[
 \|h^{(n)}\|_Y\lesssim\|h-\Pi_{n-1}h\|_{q,\ell}
                         \lesssim\|h\|_{q,\ell}.
\]
If $K:\F_m^p\to\F_\ell^q$ is compact, so is $D^nK$ as an operator
into $Y$. Thus
\[
 \|J\|_{\mathrm e}\le\|D^n(T-K)\|
                \lesssim\|T-K\|,
 \qquad \|J\|_{\mathrm e}\lesssim\|T\|_{\mathrm e}.
\]
For the reverse estimate, $TR_t$ is compact, and
\eqref{eq:integralnorm} gives
\[
 \|T\|_{\mathrm e}^q
 \le\limsup_{t\to1}\|T(\Id-R_t)\|^q
 \lesssim\limsup_{t\to1}\|J(\Id-R_t)\|^q.
\]
Combining \eqref{eq:essentialmeasure} with \eqref{eq:holomorphicapprox},
\begin{equation}\label{eq:integralessmass}
 \|T\|_{\mathrm e}^q
       \asymp\limsup_{|w|\to\infty}B_1\nu(w).
\end{equation}
The two approximations establish the comparison with compact operators in the holomorphic target. By \eqref{eq:Agdisk} and the density comparison
$d\nu\asymp A_{\mathbf g}^q\,dA$, equation~\eqref{eq:Agdisk} also gives
\[
 \limsup_{|w|\to\infty}B_1\nu(w)
             \asymp\limsup_{|z|\to\infty}A_{\mathbf g}(z)^q.
\]
Indeed, the first inequality in \eqref{eq:Agdisk} bounds the pointwise upper limit by the upper limit of disk masses, while the second bounds the latter by a constant times $\lim_{R\to\infty}\sup_{|z|>R}A_{\mathbf g}(z)^q$. Together with \eqref{eq:integralessmass}, this proves $\|T\|_{\mathrm e}\asymp\limsup_{|z|\to\infty}A_{\mathbf g}(z)$.

If $d<0$, Theorem~\ref{thm:integral} forces $\mathbf g=0$ and all assertions follow. Suppose $d\ge0$. Polynomial rigidity gives the unique decomposition
\[
 g_j(z)=c_jz^{d-j}+r_j(z),\quad \deg r_j+j<d
       \quad(j\in\mathcal I_d),
\]
and $g_j=0$ for $j\notin\mathcal I_d$. Here $c_j$ may be zero. Uniformly in $\vartheta$,
\begin{align*}
 R^{-d}P_{\mathbf g(R\ee^{i\vartheta})}
                         (R\ee^{-i\vartheta})
 &=\sum_{j\in\mathcal I_d}c_j\ee^{i(d-2j)\vartheta}+O(R^{-1}),\\
 R^{-d}P_{\mathbf g(R\ee^{i\vartheta})}^{(k)}
                         (R\ee^{-i\vartheta})
 &=O(R^{-k}),\qquad 1\le k\le N.
\end{align*}
For $d=0$, the first identity is exact and the higher derivatives
vanish. The critical symbol is
\[
 \sigma_{\mathbf g}(\ee^{i\vartheta})
   =\ee^{in\vartheta}\sum_{j\in\mathcal I_d}
                                  c_j\ee^{i(d-2j)\vartheta}.
\]
Multiplication by $\ee^{in\vartheta}$ does not change the modulus. By \eqref{eq:jet}, $S_{\mathbf g}$ is a finite sum of absolute values. Since $R^d/(1+R)^d\to1$ and every scaled derivative term of positive order tends uniformly to zero, the reverse triangle inequality, applied termwise, gives
\[
 \sup_{\vartheta}\left|
 A_{\mathbf g}(R\ee^{i\vartheta})
            -|\sigma_{\mathbf g}(\ee^{i\vartheta})|\right|\to0.
\]
It follows that
$\limsup_{|z|\to\infty}A_{\mathbf g}(z)=\max_{\mathbb T}|\sigma_{\mathbf g}|$,
proving \eqref{eq:essentialintegral}. The operator determined by the
remainders $(r_j)$ is compact by Theorem~\ref{thm:integral}.
Since $n+d-2j$ takes distinct values for distinct $j$, Fourier
orthogonality yields
\[
 \frac1{2\pi}\int_0^{2\pi}
      |\sigma_{\mathbf g}(\ee^{i\vartheta})|^2\dd\vartheta
       =\sum_{j\in\mathcal I_d}|c_j|^2.
\]
Hence the operator formed from the critical terms is compact if and only if all its coefficients $c_j$ vanish. Conversely, every choice of critical coefficients defines a bounded operator. These observations give the dimension assertion and the criterion for compact differences. In particular, 
\[
 \max_{j\in\mathcal I_d}|c_j|
 \le\max_{\mathbb T}|\sigma_{\mathbf g}|
 \le\sum_{j\in\mathcal I_d}|c_j|,
\]
where an empty index set gives zero on both sides.
\end{proof}

\Needspace{8\baselineskip}
We next compute the exact essential norm in the Hilbert case using normalized monomials.
\begin{proof}[Proof of Theorem \ref{thm:essentialintegral}: the Hilbert case]
Set $\beta_m(k)=\|z^k\|_{2,m}^2$. Radial integration and the binomial formula give the exact expression
\begin{align}
 \beta_m(k)
 &=2\int_0^\infty r^{2k+1}(1+r)^{2m}\ee^{-r^2}\dd r\notag\\
 &=\sum_{v=0}^{2m}\binom{2m}{v}\Gamma(k+v/2+1)\notag\\
 &=\Gamma(k+m+1)\big(1+O(k^{-1/2})\big)\quad(k\to\infty).
 \label{eq:monomialmass}
\end{align}
Here $\Gamma$ is the gamma function. The last line follows by applying $\Gamma(k+\alpha)/\Gamma(k+\beta) =k^{\alpha-\beta}(1+O(k^{-1}))$ to the finitely many terms. For $m=0$, the gamma formula is exact. The normalized monomials $e_k^{(m)}=z^k/\sqrt{\beta_m(k)}$ form an orthonormal basis: for $f(z)=\sum_{k\ge0}a_kz^k$, angular orthogonality and Tonelli give
\[
                         \|f\|_{2,m}^2=\sum_{k\ge0}|a_k|^2\beta_m(k).
\]
Thus $(e_k^{(m)})_{k\ge0}$ is an orthonormal basis. Let
$U_m:\ell^2(\Nzero)\to\F_m^2$ be the associated unitary.

Take a monomial $a_{h,j}z^h$ from $g_j$, and set $s=n+h-j$.
For $k\ge j$, differentiating $e_k^{(m)}$ and integrating $n$ times
gives
\begin{align}
 I^n\!\left(z^h(e_k^{(m)})^{(j)}\right)
 &=\frac{k!}{(k-j)!}
   \frac{(k-j+h)!}{(k-j+h+n)!}
       \frac{z^{k+s}}{\sqrt{\beta_m(k)}}\notag\\
 &=w_k(h,j)e_{k+s}^{(\ell)},\notag\\
 w_k(h,j)
 &=\frac{k!}{(k-j)!}
   \frac{(k-j+h)!}{(k-j+h+n)!}
       \sqrt{\frac{\beta_\ell(k+s)}{\beta_m(k)}}\notag\\
 &=k^{(h+j-d)/2}\big(1+O(k^{-1/2})\big).
 \label{eq:bandweights}
\end{align}
For $k<j$ the derivative vanishes. The asymptotic weight in the
last line follows first from the factorial ratios
\[
 \frac{k!}{(k-j)!}=k^j(1+O(k^{-1})),\qquad
 \frac{(k-j+h)!}{(k-j+h+n)!}=k^{-n}(1+O(k^{-1})).
\]
In addition, \eqref{eq:monomialmass} implies
\begin{align*}
 \frac{\beta_\ell(k+s)}{\beta_m(k)}
 &=\frac{\Gamma(k+s+\ell+1)}{\Gamma(k+m+1)}
                                  (1+O(k^{-1/2}))\\
 &=k^{s+\ell-m}(1+O(k^{-1/2})).
\end{align*}
Taking the positive square root and multiplying the three factors proves \eqref{eq:bandweights}; the exponent simplifies because
\[
 j-n+\frac{s+\ell-m}{2}=\frac{h+j-d}{2}.
\]
The factorial arguments in \eqref{eq:bandweights} are nonnegative
when $k\ge j$. If $s<0$, only finitely many initial basis vectors
require separate treatment.

Let $S$ denote the unilateral shift on $\ell^2(\Nzero)$, and define
$V_s=S^s$ for $s\ge0$ and $V_s=(S^*)^{-s}$ for $s<0$.
When $h+j<d$, \eqref{eq:bandweights} shows that the weights tend to
zero. Write $W_{h,j}$ for this weighted shift. To verify its
compactness, let $P_M$ project onto
$\Span\{e_0,\ldots,e_{M-1}\}$. For $M\ge j$, the images of distinct
basis vectors are orthogonal, so
\[
 \begin{aligned}
 \|W_{h,j}(\Id-P_M)x\|^2
   &=\sum_{k\ge M}|w_k(h,j)|^2|x_k|^2\\
   &\le\sup_{k\ge M}|w_k(h,j)|^2\|x\|^2.
 \end{aligned}
\]
Thus the finite-rank operators $W_{h,j}P_M$ converge in norm to $W_{h,j}$. If $h+j=d$, the weights tend to one. The identical argument with $w_k(h,j)-1$ proves that $W_{h,j}-V_s$ is compact, apart from a finite number of initial columns, which also give finite rank. There are only finitely many monomials. Therefore
\[
 \begin{aligned}
 U_\ell^{-1}TU_m
   &=A_\sigma+K_0,
\end{aligned}
\]
where
\[
 \begin{aligned}
 A_\sigma
   &=\sum_{j\in\mathcal I_d}c_jV_{n+d-2j},
   \qquad K_0\text{ compact}.
 \end{aligned}
\]
Identify $\ell^2(\Nzero)$ with the Hardy space $H^2(\mathbb T)$, viewed as the closed subspace of $L^2(\mathbb T,d\vartheta/(2\pi))$ generated by the nonnegative Fourier modes. Under this identification, $A_\sigma$ is the compression of
multiplication by $\sigma=\sigma_{\mathbf g}$. Hence
\[
                  \|A_\sigma\|_{\mathrm e}\le\|A_\sigma\|
                                     \le\max_{\mathbb T}|\sigma|.
\]

To obtain the reverse estimate, fix $\zeta\in\mathbb T$ and positive
integers $M,L$. In the standard basis $(e_k)$, define
\[
              x_{M,L}=L^{-1/2}\sum_{k=M}^{M+L-1}\zeta^{-k}e_k.
\]
These vectors have norm one. For every fixed integer $s$, if $M>|s|$ and $L>|s|$, cancellation on the common part of the two blocks gives
\[
 \|(V_s-\zeta^s)x_{M,L}\|^2\le\frac{2|s|}{L}.
\]
Summing over the critical indices gives
\begin{equation}\label{eq:blockestimate}
 \|(A_\sigma-\sigma(\zeta)\Id)x_{M,L}\|
 \le L^{-1/2}\sum_{j\in\mathcal I_d}|c_j|
                              \sqrt{2|n+d-2j|}.
\end{equation}
For fixed $L$, the support of $x_{M,L}$ moves to infinity as
$M\to\infty$. Hence $x_{M,L}$ converges weakly to zero in $\ell^2$.
Every compact $K$ therefore satisfies $Kx_{M,L}\to0$. We obtain
\begin{align*}
 \|A_\sigma-K\|
 &\ge\limsup_{M\to\infty}\|(A_\sigma-K)x_{M,L}\|\\
 &\ge|\sigma(\zeta)|-
       L^{-1/2}\sum_{j\in\mathcal I_d}|c_j|\sqrt{2|n+d-2j|}.
\end{align*}
Letting $L\to\infty$ and then taking the supremum over $\zeta$ and
the infimum over compact $K$, we find
$\|A_\sigma\|_{\mathrm e}=\max_{\mathbb T}|\sigma|$.
Equation~\eqref{eq:toeplitzmodel} and unitary invariance now give
\eqref{eq:essentialexact}.
\end{proof}

\begin{proof}[Proof of Theorem \ref{thm:fredholmsymbol}]
Unitary equivalence and compact perturbations preserve Fredholmness and the index \cite{Douglas}. By \eqref{eq:toeplitzmodel}, it therefore suffices to consider the Laurent Toeplitz operator $A_\sigma$, where $\sigma=\sigma_{\mathbf g}$.

Suppose that $\sigma(\zeta_0)=0$ at some point of $\mathbb T$. In the block construction above, take $\zeta=\zeta_0$ and $M=L=\nu$, with $\nu$ larger than every absolute shift index, and write $x_\nu=x_{\nu,\nu}$. These unit vectors converge weakly to zero because their supports escape to infinity. By \eqref{eq:blockestimate}, $A_\sigma x_\nu\to0$. If $A_\sigma$ were Fredholm, there would be a bounded operator $B$ and a compact operator $K_1$ such that $BA_\sigma=\Id+K_1$. It would follow that
\[
 1=\|x_\nu\|\le\|B\|\|A_\sigma x_\nu\|+\|K_1x_\nu\|
                         \longrightarrow0.
\]
This contradiction also covers $\sigma\equiv0$.

Assume now that $\sigma$ has no zero on $\mathbb T$. Choose
$r\in\Nzero$ such that $P(\zeta)=\zeta^r\sigma(\zeta)$ is a polynomial, and factor it as
\[
 P(\zeta)=c\prod_{v=1}^{M}(\zeta-a_v),\qquad c\ne0,
 \quad |a_v|\ne1.
\]
The roots are counted with multiplicity; a nonzero constant polynomial corresponds to $M=0$. The relation $S^*S=\Id$ gives, for every integer $v\ge0$,
\[
 (S^*)^rS^v=
 \begin{cases}
 S^{v-r},&v\ge r,\\
 (S^*)^{r-v},&v<r,
 \end{cases}
\]
so $(S^*)^rS^v=V_{v-r}$. Applying this identity term by term to $P$ gives the exact
factorization
\[
 A_\sigma=(S^*)^rP(S)
          =c(S^*)^r\prod_{v=1}^{M}(S-a_v\Id).
\]
When $|a_v|>1$, a Neumann series makes $S-a_v\Id$ invertible.
For $|a_v|<1$, we instead have
\[
 \|(S-a_v\Id)x\|\ge(1-|a_v|)\|x\|,
 \qquad
 \ker(S^*-\bar a_v\Id)
       =\Span\{(1,\bar a_v,\bar a_v^2,\ldots)\}.
\]
Thus $S-a_v\Id$ has closed range, zero kernel and range codimension one. Its index is $-1$. Also $\operatorname{ind}(S^*)^r=r$. The product rule for the index and the argument principle yield
\begin{align*}
 \operatorname{ind}A_\sigma
 &=r-\#\{v:|a_v|<1\}
  =-\operatorname{wind}(\zeta^{-r}P(\zeta),0)
  =-\operatorname{wind}(\sigma,0).
\end{align*}
We have proved \eqref{eq:fredholmcriterion} and
\eqref{eq:indexsymbol} for Laurent polynomial symbols, including those
with negative powers.

When $m=\ell$, the unitary maps on the domain and codomain coincide.
Consequently,
\[
 U_m^{-1}(T-\lambda\Id)U_m=A_{\sigma-\lambda}+K.
\]
Apply the preceding argument to the Laurent polynomial $\sigma-\lambda$. This proves both assertions in \eqref{eq:essspectrum}.
\end{proof}

\begin{example}\label{ex:circle}
On $\F_m^2$, let $n=N=1$ and
\[
 T f=I\big((\alpha z+\gamma)f+\beta f'\big),
                    \qquad \alpha,\beta,\gamma\in\C.
\]
Then $\sigma(\zeta)=\alpha\zeta^2+\beta$, so
\[
 \|T\|_{\mathrm e}=|\alpha|+|\beta|,
 \qquad \sigma_{\mathrm e}(T)=\beta+\alpha\mathbb T.
\]
If $\alpha\ne0$, Fredholmness is equivalent to $|\alpha|\ne|\beta|$, and
\[
 \operatorname{ind}T=
 \begin{cases}
 -2,&|\beta|<|\alpha|,\\
 0,&|\beta|>|\alpha|.
 \end{cases}
\]
If $\alpha=0$, Fredholmness holds exactly when $\beta\ne0$, with index zero. The coefficient $\gamma$ contributes only a compact term. Index zero here does not imply invertibility: every $Tf$ vanishes at the origin.
\end{example}

\begin{example}\label{ex:negative}
Different source and target Sobolev orders can produce negative
shifts. Let $n=1$, $N=4$, $m=3$, $\ell=0$, $g_4=1$, and $g_j=0$
for $j<4$. Then
\[
 Tf=I(f^{(4)})=f^{(3)}-f^{(3)}(0),\qquad
 \sigma(\zeta)=\zeta^{-3},\qquad \|T\|_{\mathrm e}=1.
\]
The kernel is the four-dimensional space of polynomials of degree at
most three. The range consists of the functions
$h\in\F_0^2$ with $h(0)=0$: Lemma~\ref{lem:LP} gives
$I^3h\in\F_3^2$ for each such $h$, and $T(I^3h)=h$.
The range therefore has codimension one, so
$\operatorname{ind}T=4-1=3$, as predicted by \eqref{eq:indexsymbol}.
\end{example}

\subsection{Linear Differential Equations}\label{subsec:odeproof}
We now prove Theorems~\ref{thm:ode} and \ref{thm:odehilbert} for equation \eqref{eq:ode}, under the coefficient condition \eqref{eq:odecoeff}.

The compact-operator fact below applies to the quasi-Banach range as well.

\begin{lemma}\label{lem:fredholm}
Suppose that $X$ is a complete complex vector space equipped with a homogeneous quasi-norm obeying a $\theta$-triangle inequality for some $0<\theta\le1$. For a compact linear operator $K:X\to X$, injectivity of $\Id+K$ implies that it is bijective and has a bounded inverse.
\end{lemma}
\begin{proof}
Write $A=\Id+K$. If $A$ were not bounded below, we could choose $\|x_j\|=1$ with $Ax_j\to0$. Compactness gives a subsequence for which $Kx_j$ converges. Since $x_j=Ax_j-Kx_j$, this subsequence converges to a unit vector in $\ker A$, contrary to injectivity. Thus $\|Ax\|\ge c\|x\|$ for some $c>0$, and consequently $\|A^jx\|\ge c^j\|x\|$. Whenever $A^jx_k$ converges, the latter bound makes $(x_k)$ Cauchy; completeness shows that the limit of $A^jx_k$ belongs to $A^jX$. Therefore $R_j=A^jX$ is closed for each $j$.

Suppose instead that $A$ is not onto. Then $R_0\supset R_1\supset\cdots$ is strictly decreasing: if $R_j=R_{j+1}$, every $x\in X$ satisfies $A^jx=A^jAy$ for some $y\in X$, and injectivity of $A^j$ gives $x=Ay$. We may therefore choose $x_j\in R_j$ such that $\|x_j\|=1$ and $\dist(x_j,R_{j+1})>1/2$. To make this choice without local convexity, take $y\in R_j\setminus R_{j+1}$, set $d=\dist(y,R_{j+1})>0$, select $v\in R_{j+1}$ with $\|y-v\|<2d$, and normalize $y-v$. If $k>j$, then $Ax_j$ and $Kx_k=Ax_k-x_k$ both belong to $R_{j+1}$. Hence
\[
 \|Kx_j-Kx_k\|=\|-x_j+(Ax_j-Kx_k)\|>\tfrac12.
\]
This separation contradicts compactness of $K$. Thus $A$ is surjective, and the lower bound makes its inverse bounded.
\end{proof}

\begin{proof}[Proof of Theorem \ref{thm:ode}]
Let $X=\F_m^p$ and $\theta=\min\{1,p\}$. Write $T=A+K$, where
\[
 Af=\sum_{j=0}^{n-1}c_j I^n(z^{n-j}f^{(j)}),\qquad
 Kf=I^n\!\left(\sum_{j=0}^{n-1}r_jf^{(j)}\right).
\]
By Theorem~\ref{thm:integral}, each term of $A$ is bounded and $K$ is compact. The $\theta$-triangle inequality for \eqref{eq:norm} gives
\[
 \|A\|^\theta\le C_{n,m,p}^\theta\sum_{j=0}^{n-1}|c_j|^\theta,
 \quad C_{n,m,p}=\max_{0\le j<n}
       \|f\mapsto I^n(z^{n-j}f^{(j)})\|_{X\to X}.
\]
With $\varepsilon_{n,m,p}=(2C_{n,m,p})^{-1}$, the coefficient assumption implies $\|A\|<1/2$. For $M>L$, we have
\[
 \left\|\sum_{k=L+1}^M(-A)^k\right\|^\theta
       \le\sum_{k=L+1}^M\|A\|^{k\theta}\longrightarrow0
       \quad(L\to\infty).
\]
Completeness now ensures that the Neumann series defines a bounded operator $B=\sum_{k=0}^\infty(-A)^k$. Its partial sums $B_M$ satisfy
\[
 (\Id+A)B_M=B_M(\Id+A)=\Id-(-A)^{M+1}.
\]
Passing to the limit proves that $B=(\Id+A)^{-1}$. This also justifies the inversion in the quasi-Banach range.

We verify injectivity of $\Id+T$ without using smallness. If $f+Tf=0$, the first $n$ derivatives at zero and the $n$th derivative identity give
\[
 f^{(j)}(0)=0\ (0\le j<n),\qquad
 f^{(n)}+\sum_{j=0}^{n-1}g_jf^{(j)}=0.
\]
Differentiating the latter equation $r$ times at zero yields
\[
 f^{(n+r)}(0)=-\sum_{j=0}^{n-1}\sum_{s=0}^{r}
        \binom{r}{s} g_j^{(s)}(0)f^{(j+r-s)}(0).
\]
All derivatives on the right have order at most $n+r-1$. Induction on $r$ shows that every Taylor coefficient of $f$ vanishes, and hence $f=0$. Therefore
\[
\Id+T=(\Id+A)\big[\Id+(\Id+A)^{-1}K\big]
\]
is invertible by Lemma~\ref{lem:fredholm}: $(\Id+A)^{-1}K$ is compact, and the factor in brackets is injective because $\Id+T$ is injective. For any entire solution $f$ of \eqref{eq:ode}, put
\[
 H=I^nG+\sum_{j=0}^{n-1}\frac{f^{(j)}(0)}{j!}z^j\in X.
\]
Let $\widetilde f=(\Id+T)^{-1}H\in X$. Since $T\widetilde f$ has a zero of order at least $n$ at the origin, $\widetilde f^{(j)}(0)=f^{(j)}(0)$ for $0\le j<n$. Differentiating $\widetilde f+T\widetilde f=H$ $n$ times yields
\[
 \widetilde f^{(n)}+\sum_{j=0}^{n-1}g_j\widetilde f^{(j)}=G.
\]
Hence $f-\widetilde f$ satisfies the homogeneous equation with zero initial values. The Taylor-coefficient induction above, which applies to all entire functions, gives $f=\widetilde f\in X$. We did not assume in advance that the given solution $f$ belongs to $X$. If all $c_j=0$, then $A=0$ and the same argument applies without a smallness condition on the lower degree coefficients.
\end{proof}

\begin{proof}[Proof of Theorem \ref{thm:odehilbert}]
For $m=\ell$, $d=n$ and $N=n-1$, formula \eqref{eq:criticalsymbol} gives
\[
 1+\sigma_{\mathbf g}(\zeta)
       =1+\sum_{j=0}^{n-1}c_j\zeta^{2(n-j)}=Q(\zeta^2).
\]
Apply \eqref{eq:essspectrum} with $\lambda=-1$. The map $\zeta\mapsto\zeta^2$ traverses $\mathbb T$ twice, so the argument principle yields
\begin{align*}
 \operatorname{ind}(\Id+T)
 &=-\operatorname{wind}(Q(\zeta^2),0)
  =-2\operatorname{wind}(Q|_{\mathbb T},0)
  =-2Z_{\mathbb D}(Q),
\end{align*}
when $Q$ has no zeros on $\mathbb T$. We have thus obtained the Fredholm criterion and \eqref{eq:odeindex}.

The induction on Taylor coefficients in the proof of Theorem~\ref{thm:ode} shows that $\Id+T$ is injective on entire functions, independently of the coefficient sizes. If $Q$ is zero-free on $\overline{\mathbb D}$, then $\Id+T$ is Fredholm of index zero and injective. Hence it is surjective and has a bounded inverse. Conversely, invertibility implies Fredholmness with index zero. The Fredholm criterion rules out zeros on $\mathbb T$, and \eqref{eq:odeindex} then rules out zeros in $\mathbb D$. This establishes (a)$\Leftrightarrow$(b).

Given initial values $a_0,\ldots,a_{n-1}$, define
\[
 P_a(z)=\sum_{j=0}^{n-1}\frac{a_j}{j!}z^j,\qquad H=I^nG+P_a.
\]
Integrating \eqref{eq:ode} $n$ times and imposing the prescribed initial data gives $(\Id+T)f=H$. The associated first-order system has a holomorphic solution on each disk centered at zero; uniqueness makes these local solutions agree on overlapping disks, producing a unique entire solution. If (a) holds, this solution is $(\Id+T)^{-1}H\in\F_m^2$. Its stated bound follows from boundedness of the inverse and $\|P_a\|_{2,m}\le C_{n,m}\sum_j|a_j|$. Thus (a) implies (c).

To prove the converse, assume (c) and fix $h\in\F_m^2$. Set $G=h^{(n)}$ and $a_j=h^{(j)}(0)$. Taylor's formula gives $I^nG=h-P_a\in\F_m^2$, so (c) yields $f\in\F_m^2$ satisfying $(\Id+T)f=h$. Surjectivity, together with the injectivity already proved, gives (a).
\end{proof}

Rouch\'e's theorem gives the sufficient condition $\max_{|w|=1}|Q(w)-1|<1$. This sufficient condition is strictly stronger than the zero-free condition in Theorem~\ref{thm:odehilbert}, as the following example shows.
\begin{example}\label{ex:ode}
For arbitrary $a,b,c\in\C$, consider
\[
 f''+(z+c)f'+\left(\tfrac14z^2+az+b\right)f=G.
\]
Here $Q(w)=1+w+w^2/4=(1+w/2)^2$, whose only zero is $-2$. Theorem~\ref{thm:odehilbert} therefore applies for every $m\in\Nzero$ and every $G$ with $I^2G\in\F_m^2$, although
\[
 \|T\|_{\mathrm e}
       =\max_{|\zeta|=1}|\zeta^2+\zeta^4/4|=\tfrac54>1.
\]
Thus solvability depends on the zero locations of $Q$, even when the essential norm is larger than one.
\end{example}

The closed disk in Theorem~\ref{thm:odehilbert} is necessary.
For the equation $f'-zf=0$ with $f(0)=1$, we have $Q(w)=1-w$
and $f(z)=\ee^{z^2/2}$. Writing $z=x+iy$, we obtain
\[
 \|f\|_{2,m}^2
 =\frac1\pi\int_{\mathbb R^2}
       (1+\sqrt{x^2+y^2})^{2m}\ee^{-2y^2}\dd x\dd y
 =\infty.
\]
Thus a zero on $\mathbb T$ can already prevent membership for zero forcing and prescribed initial data.

\subsection{Sums of Weighted Composition--Differentiation Operators}\label{subsec:compositionproof}
We use the functions $E_{w,k}$ to prove Proposition~\ref{prop:affine}, and then prove Theorems~\ref{thm:compup} and \ref{thm:compdown}. Recall the function $\mathcal B$ from \eqref{eq:B}.

\begin{lemma}\label{lem:pointB}
If $\Lu:\F_m^p\to\F_\ell^q$ is bounded, then
\begin{equation}\label{eq:Bnecessity}
                       \mathcal B(z)\lesssim\|\Lu\|,
                       \qquad z\in\C.
\end{equation}
If $\Lu$ is compact, then $\mathcal B(z)\to0$ as $|\varphi(z)|\to\infty$.
\end{lemma}
\begin{proof}
Fix $z\in\C$, let $w=\varphi(z)$, and test $\Lu$ on the functions $E_{w,k}$ from \eqref{eq:atoms}. Evaluating the derivatives of $E_{w,k}$ at $w$ gives
\[
 \Lu E_{w,k}(z)
   =\frac{\ee^{|w|^2/2}}{(1+|w|)^m}
                     P_{\mathbf u(z)}^{(k)}(\bar w).
\]
Indeed, the product rule applied to \eqref{eq:atoms} gives
\[
 E_{w,k}^{(j)}(w)=
 \begin{cases}
 0,&j<k,\\[2pt]
 \displaystyle\frac{j!}{(j-k)!}\,\bar w^{j-k}
       \frac{\ee^{|w|^2/2}}{(1+|w|)^m},&j\ge k.
 \end{cases}
\]
Summing $u_j(z)E_{w,k}^{(j)}(w)$ over $j$ gives the displayed identity. The point-evaluation estimate in $\F_\ell^q$ from Lemma~\ref{lem:local} now yields
\begin{align*}
 (1+|z|)^\ell\ee^{-|z|^2/2}|\Lu E_{w,k}(z)|
 &\lesssim\|\Lu E_{w,k}\|_{q,\ell}\\
 &\le\|\Lu\|\|E_{w,k}\|_{p,m}\\
 &\lesssim\|\Lu\|.
\end{align*}
Summing these inequalities over $k$ yields \eqref{eq:Bnecessity}. If $\Lu$ is compact and $|w|\to\infty$, Lemma~\ref{lem:compact} gives $\|\Lu E_{w,k}\|_{q,\ell}\to0$ for every $k$; the same estimates then show that $\mathcal B(z)\to0$.
\end{proof}

\begin{proof}[Proof of Proposition \ref{prop:affine}]
Assume $\mathbf u\not\equiv0$ and choose the largest $K$ with $u_K\not\equiv0$. Because $P_{\mathbf u(z)}^{(K)}(\xi)=K!u_K(z)$, Lemma~\ref{lem:pointB} yields
\begin{equation}\label{eq:affinegrowth}
 |u_K(z)|\frac{(1+|z|)^\ell}{(1+|\varphi(z)|)^m}
          \ee^{(|\varphi(z)|^2-|z|^2)/2}\le C.
\end{equation}
Away from the zeros of $u_K$, taking logarithms gives
\[
 \tfrac12|\varphi(z)|^2
 \le\tfrac12|z|^2+m\log(1+|\varphi(z)|)
       -\ell\log(1+|z|)+C-\log|u_K(z)|.
\]
To average the logarithmic term, write $u_K(z)=z^t h(z)$, where $t\in\Nzero$ and $h(0)\ne0$. Jensen's formula then implies, for $R\ge1$,
\[
 \frac1{2\pi}\int_0^{2\pi}\log|u_K(R\ee^{i\vartheta})|\dd\vartheta
 \ge t\log R+\log|h(0)|\ge\log|h(0)|.
\]
The logarithmic singularity at a zero on the circle is integrable, and the inequality at such a radius follows by continuity in $R$. Define
\[
 M_2(R)^2=\frac1{2\pi}\int_0^{2\pi}
                         |\varphi(R\ee^{it})|^2\dd t.
\]
Using $m\log(1+x)\le x^2/4+C_m$ to absorb the logarithmic growth, we find
\[
                        M_2(R)^2\le2R^2+C.
\]
Write $\varphi(z)=\sum_{j=0}^{\infty}a_jz^j$. Orthogonality of the exponentials on $|z|=R$ gives
\[
                    \sum_{j=0}^{\infty}|a_j|^2R^{2j}
                                  \le2R^2+C.
\]
For every $j\ge2$, the left side contains $|a_j|^2R^{2j}$, so letting $R\to\infty$ forces $a_j=0$. Therefore $\varphi(z)=az+b$. Take logarithms in \eqref{eq:affinegrowth} and average again on $|z|=R$. The mean of $|az+b|^2$ is $|a|^2R^2+|b|^2$; the logarithmic weight terms are $O(\log R)$, and the lower bound in Jensen's formula controls the mean of $-\log|u_K|$ from above. Consequently,
\[
          |a|^2R^2+|b|^2\le R^2+O(\log R).
\]
Dividing by $R^2$ and sending $R\to\infty$ gives $|a|\le1$.

If $a=0$, then
\[
                       \Lu f=\sum_{j=0}^Nu_j f^{(j)}(b).
\]
Each derivative evaluation $f\mapsto f^{(j)}(b)$ is bounded on the source space. Thus, if every $u_j$ belongs to $\F_\ell^q$, the displayed operator is bounded and has finite rank. For necessity, use $f_k(\zeta)=(\zeta-b)^k/k!$, $0\le k\le N$. These polynomials lie in every source space and give $\Lu f_k=u_k$, so each $u_k$ must be in the target space. If $\mathbf u\equiv0$, the claim is immediate.
\end{proof}

\begin{proof}[Proof of Theorem \ref{thm:compup}]
Lemma~\ref{lem:pointB} proves necessity in (a) and the bound $\|\mathcal B\|_\infty\lesssim\|\Lu\|$. To prove sufficiency, define, for $0\le k\le N$,
\[
 H_k(w)=(1+|w|)^m|Q_kf(w)|.
\]
Use the algebraic identity \eqref{eq:algebra} with coefficient vector $\mathbf u(z)$ and evaluation point $w=\varphi(z)$. Multiplication by the target weight gives the exact equality
\begin{align*}
 &\Lu f(z)(1+|z|)^\ell\ee^{-|z|^2/2}\\
 &\quad=\frac{(1+|z|)^\ell}{(1+|\varphi(z)|)^m}
              \ee^{(|\varphi(z)|^2-|z|^2)/2}
       \sum_{k=0}^N\frac{P_{\mathbf u(z)}^{(k)}
                                      (\overline{\varphi(z)})}{k!}
             (1+|\varphi(z)|)^m Q_kf(\varphi(z)).
\end{align*}
Taking absolute values and applying the finite-sum inequality, we obtain
\begin{equation}\label{eq:compositionmajorant}
 |\Lu f(z)|(1+|z|)^\ell\ee^{-|z|^2/2}
        \lesssim\mathcal B(z)\sum_{k=0}^N H_k(\varphi(z)).
\end{equation}
Let $M=\|\mathcal B\|_\infty<\infty$. Raise the preceding bound to the $q$th power, apply the finite-sum inequality, and change variables by $w=az+b$. Then
\begin{align*}
 \|\Lu f\|_{q,\ell}^q
 &\lesssim M^q\sum_{k=0}^N\int_{\C}H_k(az+b)^q\dd A(z)\\
 &=\frac{M^q}{|a|^2}\sum_{k=0}^N
                   \int_{\C}(1+|w|)^{mq}|Q_kf(w)|^q\dd A(w)\\
 &\lesssim\frac{M^q}{|a|^2}\|f\|_{q,m}^q
 \lesssim\frac{M^q}{|a|^2}\|f\|_{p,m}^q.
\end{align*}
The equality uses $dA(w)=|a|^2dA(z)$; the next estimate is \eqref{eq:Qglobal}, and the last follows from $p\le q$ and Lemma~\ref{lem:local}. This proves sufficiency in (a) and the stated upper norm bound.

We turn to (b). Because $a\ne0$, the inequality $|\varphi(z)|\ge|a||z|-|b|$ shows that $|\varphi(z)|\to\infty$ as $|z|\to\infty$. Lemma~\ref{lem:pointB} therefore makes $\mathcal B(z)\to0$ necessary for compactness. Conversely, suppose this limit holds. Continuity of $\mathcal B$ implies boundedness, so part (a) applies. Split the integral defining $\|\Lu f\|_{q,\ell}^q$ into $|z|\le R$ and $|z|>R$. On the first region the coefficients are bounded and $|\varphi(z)|\le R_a:=|a|R+|b|$; on the second, apply \eqref{eq:compositionmajorant} and the same change of variables as above. For $f\in\F_m^p$ this yields
\begin{align}
 \|\Lu f\|_{q,\ell}^q
 &\le C_R\sum_{k=0}^N\sup_{|w|\le R_a}|f^{(k)}(w)|^q\notag\\
 &\quad+C|a|^{-2}\sup_{|z|>R}\mathcal B(z)^q\|f\|_{p,m}^q.
 \label{eq:compsplit}
\end{align}
Here $C_R$ is finite for each fixed disk, whereas $C$ is independent of $R$ and $f$. If $(f_j)$ is bounded in $\F_m^p$ and converges uniformly to zero on compact sets, Cauchy's formula implies that all derivatives appearing in the first line converge uniformly to zero on $\overline{D(0,R_a)}$. Hence, for fixed $R$,
\[
 \limsup_{j\to\infty}\|\Lu f_j\|_{q,\ell}^q
 \lesssim |a|^{-2}\sup_{|z|>R}\mathcal B(z)^q
                         \sup_j\|f_j\|_{p,m}^q.
\]
The assumed decay of $\mathcal B$ makes the right side tend to zero as $R\to\infty$. Lemma~\ref{lem:compact} now proves compactness.

For the essential norm, assume $L=\Lu$ is bounded and write
\[
 M_\infty=\limsup_{|z|\to\infty}\mathcal B(z).
\]
Let $K:\F_m^p\to\F_\ell^q$ be compact. Taking $w=\varphi(z)$ in the evaluation identity from Lemma~\ref{lem:pointB}, we obtain
\begin{align*}
 \mathcal B(z)
 &=\sum_{k=0}^N|LE_{w,k}(z)|(1+|z|)^\ell\ee^{-|z|^2/2}\\
 &\le\sum_{k=0}^N|(L-K)E_{w,k}(z)|
                                  (1+|z|)^\ell\ee^{-|z|^2/2}\\
 &\quad+\sum_{k=0}^N|KE_{w,k}(z)|
                                  (1+|z|)^\ell\ee^{-|z|^2/2}\\
 &\lesssim\|L-K\|+
       \sup_{\|f\|_{p,m}\le C_0}
                         |Kf(z)|(1+|z|)^\ell\ee^{-|z|^2/2},
\end{align*}
where $C_0=\max_{0\le k\le N}\sup_w\|E_{w,k}\|_{p,m}<\infty$.
By \eqref{eq:compactpointtails}, the last term vanishes uniformly as $|z|\to\infty$. Taking the upper limit and then the infimum over compact operators $K$, we obtain
\[
 M_\infty\lesssim\|L\|_{\mathrm e}.
\]

To obtain the reverse estimate, note that $LR_t$ is compact by Lemma~\ref{lem:dilations}. Apply \eqref{eq:compsplit} to $h=(\Id-R_t)f$ for $\|f\|_{p,m}\le1$ and $t\ge1/2$. The quantities $\|h\|_{p,m}$ remain uniformly bounded, while \eqref{eq:dilationlocal} makes the first line of \eqref{eq:compsplit} tend uniformly to zero for each fixed $R$. Therefore
\[
 \limsup_{t\to1}\|L(\Id-R_t)\|^q
       \lesssim |a|^{-2}\sup_{|z|>R}\mathcal B(z)^q.
\]
Since $L-LR_t=L(\Id-R_t)$ and $LR_t$ is compact, sending $R\to\infty$ gives
\[
 \|L\|_{\mathrm e}^q
 \le\limsup_{t\to1}\|L(\Id-R_t)\|^q
       \lesssim |a|^{-2}M_\infty^q.
\]
Together with the lower bound, this is the claimed essential-norm estimate.
\end{proof}

\begin{proof}[Proof of Theorem \ref{thm:compdown}]
We first prove the necessity of (c). Assume first that $L=\Lu$ is bounded. Set $\psi(w)=(w-b)/a$ and $v_j(w)=u_j(\psi(w))$. The substitution $w=az+b$ in the target norm yields
\begin{align}
 \|Lf\|_{q,\ell}^q
 &=\frac q{2\pi}\int_{\C}
       \left|\sum_{j=0}^Nu_j(z)f^{(j)}(az+b)\right|^q
                   (1+|z|)^{\ell q}\ee^{-q|z|^2/2}\dd A(z)\notag\\
 &=\frac q{2\pi|a|^2}\int_{\C}
       \left|\sum_{j=0}^Nv_j(w)f^{(j)}(w)\right|^q
       (1+|\psi(w)|)^{\ell q}\ee^{-q|\psi(w)|^2/2}\dd A(w).
 \label{eq:pullback}
\end{align}
Define
\[
 d\mu_\varphi(w)=|a|^{-2}(1+|\psi(w)|)^{\ell q}
                \ee^{q(|w|^2-|\psi(w)|^2)/2}\dd A(w).
\]
The resulting measure is positive and locally finite. With this notation, \eqref{eq:pullback} becomes
\[
 \|Lf\|_{q,\ell}^q
 =\frac q{2\pi}\int_{\C}
       \left|\sum_{j=0}^Nv_j(w)f^{(j)}(w)\right|^q
                        \ee^{-q|w|^2/2}\dd\mu_\varphi(w).
\]
Consequently $J_{\mathbf v,\mu_\varphi}$ is bounded, and the preceding identity gives
\[
 \|L\|=\left(\frac{q}{2\pi}\right)^{1/q}
           \|J_{\mathbf v,\mu_\varphi}\|.
\]
The associated measure $\nu_\varphi$ has density
\begin{align}
 d\nu_\varphi(w)
 &=|a|^{-2}\frac{(1+|\psi(w)|)^{\ell q}}{(1+|w|)^{mq}}
       \ee^{q(|w|^2-|\psi(w)|^2)/2}
       \sum_{k=0}^N|P_{\mathbf v(w)}^{(k)}(\bar w)|^q\dd A(w)\notag\\
 &\asymp |a|^{-2}\mathcal B(\psi(w))^q\dd A(w).
 \label{eq:pullbackmeasure}
\end{align}
The superscript $(k)$ here means differentiation in the polynomial variable. In particular, it does not differentiate $u_j\circ\psi$. Put $r=\frac{p}{p-q}$; then $s=qr$, and Theorem~\ref{thm:carlesondown} gives
\begin{equation}\label{eq:pullbackaverage}
 B_1\nu_\varphi\in L^r,
 \qquad \|B_1\nu_\varphi\|_{L^r}\lesssim\|L\|^q.
\end{equation}

We next derive the integrability condition in (c) from these disk averages. Set $A_\varphi(w)=\mathcal B(\psi(w))$. Because $\varphi(\psi(w))=w$, we can write
\[
 \begin{aligned}
 A_\varphi(w)
   &=W_\varphi(w)\ee^{\beta_\varphi(w)}
      \sum_{k=0}^N|P_{\mathbf v(w)}^{(k)}(\bar w)|,
 \qquad W_\varphi(w)
   &=\frac{(1+|\psi(w)|)^\ell}{(1+|w|)^m},
 \end{aligned}
\]
where
\[
 \beta_\varphi(w)
 =\frac{1-|a|^{-2}}2|w|^2
       +|a|^{-2}\operatorname{Re}(\bar b w)
       -\frac{|b|^2}{2|a|^2}.
\]
When $|w-\zeta|<1$, we have $|\psi(w)-\psi(\zeta)|<|a|^{-1}$ and $1+|w|\asymp1+|\zeta|$. Thus $W_\varphi$ varies by at most a fixed factor on these disks. Apply Lemma~\ref{lem:submeanjet} with $\eta(w)=w$ and exponent $q$ to obtain
\[
 A_\varphi(w)^q
   \lesssim_\varphi\int_{D(w,1)}A_\varphi(\zeta)^q\dd A(\zeta)
   \asymp |a|^2 B_1\nu_\varphi(w).
\]
The final comparison follows by integrating \eqref{eq:pullbackmeasure} over $D(w,1)$. Raise the pointwise inequality to the power $r$ and integrate. Using \eqref{eq:pullbackaverage} and $w=az+b$, we find
\begin{align*}
 \int_{\C}\mathcal B(z)^s\dd A(z)
 &=|a|^{-2}\int_{\C}A_\varphi(w)^{qr}\dd A(w)\\
 &\lesssim_\varphi |a|^{2r-2}
                       \int_{\C}(B_1\nu_\varphi(w))^r\dd A(w)\\
 &\lesssim_\varphi\|L\|^{qr}=\|L\|^s.
\end{align*}
We have proved (a)$\Rightarrow$(c) and the lower bound for the operator norm.

Now assume $\mathcal B\in L^s$. For $f\in\F_m^p$, define
\[
 H_f(w)=\sum_{k=0}^N(1+|w|)^m|Q_kf(w)|.
\]
The elementary inequality for a finite sum, followed by \eqref{eq:Qglobal}, gives
\[
 \int_{\C}H_f(w)^p\dd A(w)\lesssim\|f\|_{p,m}^p.
\]
Since $q/s+q/p=1$, the exponents $s/q$ and $p/q$ are conjugate and exceed one. Apply H\"older's inequality to \eqref{eq:compositionmajorant}, then use the Jacobian $|a|^{-2}$, to get
\begin{align*}
 \|Lf\|_{q,\ell}^q
 &\lesssim\int_{\C}\mathcal B(z)^qH_f(az+b)^q\dd A(z)\\
 &\le\left(\int_{\C}\mathcal B(z)^s\dd A(z)\right)^{q/s}
       \left(\int_{\C}H_f(az+b)^p\dd A(z)\right)^{q/p}\\
 &=|a|^{-2q/p}\|\mathcal B\|_{L^s}^q
                \left(\int_{\C}H_f(w)^p\dd A(w)\right)^{q/p}\\
 &\lesssim |a|^{-2q/p}\|\mathcal B\|_{L^s}^q\|f\|_{p,m}^q.
\end{align*}
This gives (c)$\Rightarrow$(a) and the upper norm bound. The argument uses H\"older's inequality with exponents greater than one and therefore also covers $p<1$ or $q<1$.

It remains to prove compactness. Let $(f_j)$ be bounded in $\F_m^p$ and converge to zero uniformly on compact subsets of $\C$. By Cauchy's formula, each $f_j^{(k)}$ tends uniformly to zero on the compact set $\varphi(\overline{D(0,R)})$; thus the integral over $|z|\le R$ tends to zero. On $|z|>R$, the preceding H\"older argument gives
\begin{align*}
 &\int_{|z|>R}|Lf_j(z)|^q(1+|z|)^{\ell q}
                                      \ee^{-q|z|^2/2}\dd A(z)\\
 &\quad\lesssim |a|^{-2q/p}
       \left(\int_{|z|>R}\mathcal B(z)^s\dd A(z)\right)^{q/s}
                                            \|f_j\|_{p,m}^q.
\end{align*}
Because $\mathcal B\in L^s$, the tail integral tends to zero as $R\to\infty$, uniformly for the bounded sequence $(f_j)$. First choose $R$ to control this tail, and then let $j\to\infty$ on the fixed disk. We obtain $\|Lf_j\|_{q,\ell}\to0$, so Lemma~\ref{lem:compact} proves (c)$\Rightarrow$(b). Since (b)$\Rightarrow$(a) is immediate, the three conditions are equivalent. The preceding estimates give the norm comparison, while compactness implies $\|L\|_{\mathrm e}=0$.
\end{proof}

\begin{proof}[Proof of Corollary \ref{cor:commute}]
For constant $\varphi$, we have $D^nC_\varphi=0$, whereas $C_\varphi D^n$ has rank one. Boundedness of derivative evaluations shows that both operators are bounded and compact. Let $\varphi$ now be nonconstant. The chain rule gives
\[
 D^nC_\varphi f
       =(\varphi')^n f^{(n)}\circ\varphi
                    +\sum_{j=1}^{n-1}u_j f^{(j)}\circ\varphi.
\]
Here each $u_j$ is entire and $(\varphi')^n\not\equiv0$. If either operator is bounded, Proposition~\ref{prop:affine} implies $\varphi(z)=az+b$ with $0<|a|\le1$. For this affine map, all terms involving derivatives $f^{(j)}$ with $j<n$ vanish, leaving
\[
                        D^nC_\varphi=a^n C_\varphi D^n.
\]
This identity gives the asserted equivalences for boundedness and compactness. In the bounded case it also yields
\[
 \|D^nC_\varphi\|_{\mathrm e}
               =|a|^n\|C_\varphi D^n\|_{\mathrm e}.\qedhere
\]
\end{proof}

The preceding criteria make the corollary explicit. For $C_\varphi D^n$, the coefficient polynomial is $P(\xi)=\xi^n$, and therefore
\[
 \mathcal B(z)\asymp
 (1+|z|)^\ell(1+|\varphi(z)|)^{n-m}
                 \ee^{(|\varphi(z)|^2-|z|^2)/2}.
\]
When $0<|a|<1$, the negative quadratic part of the exponent dominates the linear and polynomial factors, making both operators compact for all $p,q,m,\ell$. For $|a|=1$ and $b\ne0$, the exponential grows on an open angular sector, so neither operator is bounded. In the remaining case $|a|=1$ and $b=0$, the conditions become
\[
 \renewcommand{\arraystretch}{1.35}
 \setlength{\arraycolsep}{9pt}
 \begin{array}{@{}lcc@{}}
 \toprule
 &\text{boundedness}&\text{compactness}\\
 \midrule
 p\le q&m-\ell\ge n&m-\ell>n\\
 q<p&m-\ell>n+2\big(\frac1q-\frac1p\big)&m-\ell>n+2\big(\frac1q-\frac1p\big).\\
 \bottomrule
 \end{array}
\]
In particular, $D^n:\F_{\ell+n}^p\to\F_\ell^p$ is bounded but not compact; increasing $m-\ell$ makes it compact. For $q<p$, equality at $m-\ell=n+2\big(\frac1q-\frac1p\big)$ gives the logarithmic divergence in \eqref{eq:radialthreshold}.

\bigskip
\noindent\textbf{Statements and Declarations}
\par\medskip
\noindent\textbf{Funding.} This work was supported by the National Natural Science Foundation of China (Grant No.~12601234).

\section*{Data availability}
No data were generated or analyzed for the research described in this article.

\section*{Competing interests}
The authors declare that they have no competing interests.

\Needspace{7\baselineskip}
\bigskip
\noindent\textbf{Jianxiang Dong and Shibo Hu}
\par\noindent$^{1}$School of Mathematics and Statistics, Tianshui Normal University,\\
Tianshui 741000, P.R. China
\par\smallskip
\noindent E-mail: \href{mailto:jianxd@tsnu.edu.cn}{jianxd@tsnu.edu.cn} (Jianxiang Dong)
\par\noindent E-mail: \href{mailto:GL_agoin@outlook.com}{\nolinkurl{GL_agoin@outlook.com}} (Shibo Hu)
\par\smallskip
\noindent$^{*}$Corresponding author: Shibo Hu.
\end{document}